\documentstyle[12pt,amsfonts,amscd]{amsart}
\newtheorem{Theorem}{Theorem}[section]
\newtheorem{Definition}[Theorem]{Definition}
\newtheorem{Proposition}[Theorem]{Proposition}

\newtheorem{Lemma}[Theorem]{Lemma}
\newtheorem{Corollary}[Theorem]{Corollary}
\theoremstyle{remark}
\newtheorem{Example}[Theorem]{Example}

\def\il{\int\limits_}

\def\eps{\varepsilon}

\def\ovr{\overline}

\def\om{\omega}
\def\Om{\Omega}
\def\al{\alpha}
\def\gm{\gamma}
\def\th{\theta}

\def\Gm{\Gamma}

\def\dl{\delta}

\def\bd{\partial}
\def\lm{\lambda}

\def\si{\sigma}

\def\sm{\setminus}
\def\sbs{\subset}

\def\wtl{\widetilde}

\def\supp{\operatorname{supp}}

\def\dist{\operatorname{dist}}

\def\reg{\operatorname{reg}}
\def\sing{\operatorname{sing}}

\def\re{{\mathbf {Re\,}}}
\def\im{\mathbf {Im\,}}
\def\be{\begin{enumerate}}
\def\ee{\end{enumerate}}
\def\bT{\begin{Theorem}}
\def\eT{\end{Theorem}}
\def\bP{\begin{Proposition}}
\def\eP{\end{Proposition}}
\def\bD{\begin{Definition}}
\def\eD{\end{Definition}}
\def\bE{\begin{Example}}
\def\eE{\end{Example}}
\def\bL{\begin{Lemma}}
\def\eL{\end{Lemma}}
\def\bC{\begin{Corollary}}
\def\eC{\end{Corollary}}

\def\H{{\mathcal H}}

\def\aC{{\Bbb C}}

\def\dD{{\Bbb D}}
\def\rR{{\Bbb R}}

\def\E{{\mathcal E}}

\def\D{{\mathcal D}}

\begin{document}
\title{Hardy and Bergman spaces on hyperconvex domains
and their composition operators}
\author{Evgeny A. Poletsky  and Michael I. Stessin}
\address{ Department of Mathematics,  215 Carnegie Hall,
Syracuse University,  Syracuse, NY 13244, eapolets@@syr.edu}
\address{Department of Mathematics and Statistics, University at Albany,
Albany, NY 12222, stessin@@csc.albany.edu}
%\begin{abstract}
%\end{abstract}
\maketitle
\section{Introduction}
\par The theory of Hardy and Bergman spaces of analytic functions on
the unit disk is one of the most developed and useful branch of
function theory. In several variables such spaces were studied on
the unit ball (\cite{Ru1}), strongly pseudoconvex domains
(\cite{St}, \cite{BFG}) and polydisks (\cite{Ru2}). However, in
several variables none of the listed classes of domains can serve
as a ``model'' domain similar to the unit disk in one variable.

\par The goal of this paper is to develop a technique leading to a
meaningful and uniform function theory of spaces of analytic functions on a
broad class of domains in $\aC^n$. We introduce Hardy and weighted Bergman
spaces on hyperconvex domains and prove their basic properties. Our prime
focus is on geometric aspects rather than reproducing kernels and duality. The
geometry of the domain is hidden in the exhaustion functions, but  the Nevanlinna counting
functions determined by the exhaustions reveals the geometric nature
of the norm and, therefore, spaces. Since the introduced spaces are
expected to behave under holomorphic transformations of domains similar to
the classical ones, in the last part of the paper we prove some estimates
for composition operators induced by holomorphic mappings. These estimates
might be viewed as generalized embedding results. They lead to embedding
theorems for holomorphic isomorphisms.

\par Let us briefly present the content of the paper. In the classical theory,
methods of potential theory play an important role. In several
variables, where most of other methods either disappear or become
technically difficult, we still can rely on methods of
pluripotential theory (see \cite{Kl} and \cite{D3}). In this
theory the role of the Laplacian is played by the Monge--Ampere
operator. For this operator Demailly in (\cite{D1}) proved on
hyperconvex domains a fundamental Lelong-Jensen formula, which may
be viewed as an analogue of the classical Littlewood-Paley
identity. Recall that a domain $D\subset \aC^n$ is called {\it
hyperconvex}, if there is a continuous negative plurisubharmonic
function $u$ on $D$, called an {\it exhaustion function}, such
that $\lim_{z\to\bd D}u(z)=0$. To make our presentation
self-contained we stated Lelong-Jensen formula along with other
background materials in Section \ref{S:br}.

\par Lelong-Jensen formula suggests hyperconvex domains as a natural
class of domains where one might have a rich function theory. This class is
very wide: by the theorem of Demailly (\cite{D1}) every bounded
pseudoconvex domain with Lipschitz boundary is hyperconvex. At the same
time as we show the presence of Lelong-Jensen formula suffices to prove
many interesting results.

\par Using this formula in Section \ref{S:spf} we introduce Hardy and Bergman
norms on the cone of non-negative plurisubharmonic functions on a
hyperconvex domain $D\subset \aC^n$. In general, such a norm depends on the
choice of an exhausting function $u$, but we prove that if two exhausting
functions have the same rate of decay near the boundary of the domain, then
the norms are equivalent. Moreover the smallest norms are obtained when $u$
belongs to the compact class ${\mathcal E}_0(D)$, i.e., $(dd^cu)^n$ has a
compact support in $D$.  Pluricomplex Green functions are in this class.

\par In Section \ref{S:HBs} we define Hardy and weighted Bergman norms of
a holomorphic function $f$ as the norms from Section \ref{S:spf} of
$|f|^p$. Both norms depend on the choice of the exhaustion function $u$
and, consequently, Hardy and Bergman spaces depend on $u$.
%This fact makes the
%theory richer, however we do not exploit this in the paper.
The largest spaces are obtained when $u\in{\mathcal E}_0(D)$. We
show that the latter spaces coincide with spaces that were studied
on the unit ball (\cite{Ru1}), strongly pseudoconvex domains
(\cite{St}) and polydisks (in the case of Hardy spaces)
(\cite{Ru2}).

\par The choice of a pluricomplex Green function as an exhausting
function not only frequently leads to a simplification of
estimates and resulting formulas but also gives bounded point
evaluations in all cases. Using this we prove that the introduced
spaces are Banach. In the same section we also prove a version of
the Littlewood's subordination principle for holomorphic mappings
of a hyperconvex domain into the unit disk.

\par The formulas for the norms are more complicated comparatively to
the plane case and this might be viewed as an obstruction to the
new theory. Fortunately, Lelong-Jensen formula provides a
splitting for $\|f\|^p_{A^p_{u,\al}}$, namely,
$$\|f\|^p_{A^p_{u,\al}}=\il D\si_\al(u)|f|^p(dd^cu)^n+
\il{D}\gm_\al(u)dd^c|f|^p\wedge(dd^cu)^{n-1},$$ where $\si_\al$
and $\gm_\al$ are some auxiliary functions (see Section
\ref{S:ncf}). The first integral when, for example, $u\in\E_0(D)$,
is over a compact part of $D$ and can be easily estimated.

\par We define and use  the Nevanlinna counting function $N_{\al,f}(w)$
of a holomorphic function $f$ on $D$ to understand the geometry of
the second integral. Earlier, Nevanlinna counting functions in
several variables were constructed and used in  multidimensional
Nevanlinna theory by Griffiths and King \cite{GK}. Our
construction differs from Griffiths and King's in two aspects.
First, we consider the hyperbolic case instead of parabolic. The
second and the main difference is that we construct Nevanlinna
functions for an arbitrary current of a type described below,
while in \cite{GK} such functions are defined only for {\it
special exhausting } functions.

\par The geometric meaning of the Nevanlinna counting functions is similar to
the classical one. If in the classical theory they, basically, count the
number of points where $f(z)=w$ in several variables they count areas of
the sets $\{f=w\}$ in the pseudo-metric $dd^cu$.

\par To define the Nevanlinna counting functions in Section \ref{S:cvf} we
consider the current $T=u_0dd^cu_1\wedge ...\wedge dd^cu_{n-1}$, where
$u_0,...,u_{n-1}$ are continuous plurisubharmonic functions on $D$
satisfying some mild conditions. Given a compactly supported
$C^\infty$-function $\phi $ in $D$, a holomorphic function $f$ in $D$ and
$w\in \aC$ we define
$$
N_{f,T}(w,\phi)=N(w,\phi)=\int_D \phi dd^c\log|f-w|\wedge T.
$$
\par The main result of this most technical section is Theorem \ref{T:cvi},
which establishes a change of variables formula. It states that for every
holomorphic function $f$ mapping $D$ into $\Omega\subset \aC $, every
$C^\infty$-function $\phi $ compactly supported in $D$, a
$C^\infty$-function $\psi$ compactly supported in $\Omega$ and every
subharmonic function $v$ on $\Omega$, whose infinite locus does not
intersect the image of infinite locus of $u_j, \ j=0,...,n-1$, we have
$$
\int_D \phi \psi \circ f dd^c(v\circ f)\wedge T=\int_\Omega \psi
N_{f,T}(w,\phi\psi\circ f)dd^cv.
$$

\par In Section \ref{S:ncf} we set the functions $u_0,...,u_{n-1}$ equal to
an exhaustion function $u$ on $D$ and define the Nevanlinna counting
functions of order $\al$ as
$$N_{f,\al}(w)=\il{D}\gm_\al(u)(dd^cu)^{n-1}\wedge dd^c\log|f-w|.$$
The change of variables formula implies that
$$\|f\|^p_{A^p_{u,\al}}=\il D\si_\al(u)|f|^p(dd^cu)^n+
\il{\Bbb C}N_{\al,f}(w)dd^c|w|^p.$$ The latter formula generalizes
classical Littlewood-Paley identity. When a pluricomplex Green
function is chosen as the exhaustion, the formula takes exactly
the form of this identity.

\par Section \ref{S:pncf} is devoted to properties of counting functions.
First, we show their equivalence when exhausting functions are
equivalent. We also prove that, like in the one-dimensional case,
the Nevanlinna counting function satisfies Shapiro's mean value
inequality. Finally, we prove that for every holomorphic function
which takes $D$ into the unit disk ${\Bbb D}$, the classical
logarithmic estimate holds for the Nevanlinna counting function,
i.e., $N_{u,f}(w)\le c\log|w|$  for every $w\in{\Bbb D}$
sufficiently close to the boundary.

\par It is well known that the projection $(z_1,z_2)\to z_1$ induces an
isometry of the classical Bergman space in the unit disk into the Hardy
space in the unit ball in $\aC^2$. This and many other embedding theorems
for spaces of analytic functions have numerable applications to operator
theory. In the rest of the paper we investigate more general problem: how
our spaces on domains $D_1\sbs{\Bbb C}^n$ and $D_2\sbs{\Bbb C}^m$ are
transformed by a holomorphic mapping $F:D_1\to D_2$, which induces the
operator $C_Ff=f\circ F$ called the {\it composition operator.}

\par The case when $D_1=D_2={\Bbb D}$ has been
intensively investigated since 1960-s and is well understood. A
good exposition of main results of one-dimensional theory can be
found in monographs (\cite{Sha1}) and (\cite{CoM}) and references
there.

\par Contrary to this,  results in the multivariable case are sparse. It
was understood quite a while ago that the situation in several
variables is considerably harder than in the classical setting.
For instance, Littlewood subordination principle implies that
every holomorphic self-mapping of the unit disk induces a
composition operator which acts boundedely from every weighted
Bergman space $A^p_\alpha (\dD)$ into itself ($p>0, \ \alpha \geq
-1$, and, as usually, if $\alpha =-1$, the corresponding space is
the Hardy space). For domains in $\aC ^n$ this is not the case.
Even a quadratic polynomial self-mapping of the unit ball in $\aC
^n$ need not to induce a bounded operator acting on the Hardy
space. Counter-examples were constructed by Shapiro, Cima and
Wogen (\cite{CW2}) and others.

\par In the classical case conditions for a composition operator to act
continuously or compactly from $A^p_\al({\Bbb D} )$ to
$A^p_\beta({\Bbb D})$ are naturally expressed in term of
the Nevanlinna counting functions. This was discovered by Shapiro in
his paper \cite{Sha2} when $\al=\beta=-1$ and then by Smith
(\cite{Sm}) for any $\al$ and $\beta$. Roughly speaking, the
function $N_{F,\beta}(w)$ must decay as $\gm_\al(w)$ at the
boundary for the continuity of $C_F$ and faster than that for
compactness.

\par In section 8 we prove two theorems which respectively give sufficient
and necessary conditions of boundedness and compactness of a
composition operator induced by a holomorphic mapping of arbitrary
hyperconvex domains as an operator acting from one weighted
Bergman space into another. These theorems might be considered as
multidimensional analogs of the results by Shapiro and  Smith,
though our necessary and sufficient conditions in general case
seem to be different.

\par However, in the case when $D_2={\Bbb D}$ as we show in Section
\ref{S:coimud} the gap between our necessary and sufficient conditions
disappears and they become identical in the form to the conditions
of Shapiro and Smith.

\par In section \ref{S:comsp} we prove that if the domain $D_2$ is strongly
pseudoconvex then $C_F$ maps $A^p_{\rho ,\alpha}(D_2)$
continuously into $A^p_{u,n+\alpha -1}(D_1), \ \alpha \geq -1$. In
such a generality this result cannot be improved. Previously,
MacCluer and Mercer (\cite{MM}) proved this  for self-mappings of
a strongly convex domain in $\aC ^n$ and $\al=-1$.  Cima and
Mercer (\cite{CM}) extended this result to Bergman spaces on a
strongly convex domain and any $\al\ge-1$. A similar result for
mappings between balls of arbitrary dimensions (not necessarily
equal) was obtained recently by Koo and Smith (\cite{KS}).
Mappings of polydisks were considered in (\cite{SZ1}, \cite{SZ2}).
\par The proof exploits the fact that boundedness and compactness of
composition operators are naturally expressed in terms of Carleson
measures. This approach was first explicitly stated  by MacCluer
(\cite{Mc}) (earlier Carleson measures were used by Cima and Wogen
(\cite{CW1}) to give a compactness criterion for Toeplitz operators).
\par %In section \ref{S:hmcs}
The  general result stated above might be significantly improved under some
additional assumptions. We describe two such situations. The first is  when
$D_1$ is strongly pseudoconvex, $D_2={\Bbb D}$ and $F$ is holomorphic in a
neighborhood of $D_1$. In this case we use our sufficient conditions from
Section \ref{S:coimud} to show that $C_F$ acts boundedly  from $A^p_\alpha
({\Bbb D})$ into $A^p_{u,\beta}(D)$ when $\al\le\beta+(n-1)/2$ and
compactly when $\al<\beta+(n-1)/2$. An example in this section of a
quadratic polynomial $F$ shows that the result cannot be improved.

\par The second is when $F$ is a  proper holomorphic
mappings of a hyperconvex domains $D_1$ into a  hyperconvex domain $D_2$ of
the same dimension, then for every pair of exhausting function $u_1\in
{\mathcal E_0}(D_1)$ and $u_2\in {\mathcal E_0}(D_2)$ the composition
operator $C_F$ acts boundedly from $A^p_\alpha (D_2)$ into $A^p_\alpha
(D_1)$ and has a left inverse with the same properties. In particular,
every automorphism of a strongly pseudoconvex domain induces a bounded
composition operators acting on Hardy and Bergman spaces.
\section{Background results}\label{S:br}
\subsection{Differential forms and currents}
\par Let $D$ be a domain in ${\Bbb C}^n$ and let $C^\infty_0(D)$
be the space of all smooth functions on $D$ with compact supports.
A sequence $\{\phi_j\}\sbs C^\infty_0(D)$ converges to 0 if the
supports of all $\phi_j$ belong to a compact set $K\sbs D$ and the
functions $\phi_j$ with all derivatives converge uniformly to 0.
\par We denote by $\D^{p,q}(D)$ the space of all differential form
$$\om=\sum_{|I|=p,|J|=q}\om_{IJ}dz_I\wedge d\ovr z_J$$
of bidegree $(p,q)$, where $I=\{i_1,\dots,i_p\}$ and
$J=\{j_1,\dots,j_q\}$ are subsets of $\{1,\dots,n\}$,
$dz_I=dz_{i_1}\wedge\cdots\wedge dz_{i_p}$, $d\ovr z_J=d\ovr
z_{j_1}\wedge\cdots\wedge d\ovr z_{j_q}$ and $\om_{IJ}\in
C^\infty_0(D)$. Equipped with the topology of uniform convergence
on compacta with all derivatives, $\D^{p,q}$ has a structure of a
linear topological space.
\par The space $\D'_{p,q}(D)$ of continuous linear functionals on
$\D^{p,q}(D)$ is called the space of currents of bidimension
$(p,q)$ or of bidegree $(n-p,n-q)$. If $\phi\in\D'_{p,q}(D)$ then
$$\phi=\sum_{|I|=n-p,|J|=n-q}\phi_{IJ}dz_I\wedge d\ovr z_J,$$
where $\phi_{IJ}$ are distributions and the pairing $\langle
\phi,\omega \rangle $ is given by
$$
\langle \phi,\omega \rangle=\sum_{|I|=n-p,|J|=n-q} \langle
\phi_{I,J},\om_{I,J}\rangle.
$$
\par A current $\phi\in\D'_{p,p}$ is positive if
$\langle\phi,\om\rangle\ge0$ for every test form
$$\om=i\om_1\wedge\ovr\om_1\wedge\cdots\wedge
i\om_p\wedge\ovr\om_p,\qquad\om_j\in\D^{1,0}(D).$$ In this case
the coefficients $\phi_{IJ}$ are positive measures.
\par The differential of $\omega $ is defined by
$d\om=\bd\om+\ovr\bd\om$, where
$$\bd\om=
\sum\frac{\bd\om_{IJ}}{\bd z_k}dz_k\wedge dz_I\wedge d\ovr
z_J,\qquad\ovr\bd\om=\sum\frac{\bd\om_{IJ}}{\bd\ovr z_k}d\ovr
z_k\wedge dz_I\wedge d\ovr z_J.$$ The operator $d^c$ is defined by
$d^c=i(\ovr\bd-\bd)$. For $\phi \in C^2 (D)$ we have
 $$dd^c\phi=2i\sum\frac{\bd^2\phi}{\bd z_i\bd\ovr
z_j}dz_i\wedge d\ovr z_j.$$
\par Given a current $T$ we define $dT$ by the formula: $\langle
dT,\om\rangle=\langle T,d\om\rangle$ and $dd^cT$ by the formula:
$\langle dd^cT,\om\rangle=\langle T,dd^c\om\rangle$. A current $T$
is closed if $dT=0$.

Every plurisubharmonic function generates a closed positive (1,1)-current.
The following result can be found in \cite[Prop. 3.3.5]{Kl}.

\bT\label{T:pos} If $u$ is a plurisubharmonic function on $\Omega $,
then $dd^cu$ is a closed, positive (1,1)-current with measure coefficients.
\eT

The wedge product of a $(p,q)$-form $\omega$ and $(n-r,n-s)$-current
$T$ is a $(n-r+p,n-s+q)$-current $T\wedge \omega $ defined by
$$
\langle T\wedge \omega, \psi \rangle =\langle T, \psi \wedge \omega  \rangle.
$$
In some cases it is possible to make sense of the wedge product
even when $\omega $ is not smooth (and is considered as a
current). We will use three results of this kind. To formulate
them, we introduce the infinite locus $L(u)$ of a plurisubharmonic
function $u$ on $D$ as the set of points $z\in D$ such that $u$ is
not bounded on any neighborhood of $z$.
\par The first result is due to Demailly \cite[Cor. 2.3]{D2}.
\bT\label{T:dem} Suppose that $u$ is a plurisubharmonic function,
the set $L(u)$ is compact and $T$ is a closed positive current of
bidegree $(1,1)$, then the current $(dd^cu)^{n-1}\wedge T$ is well
defined and has locally finite mass on $D$. Moreover, if $\{u_k\}$
is a sequence of decreasing plurisubharmonic functions converging
to $u$, then the currents $(dd^cu_k)^{n-1}\wedge T$ converge
weak-$*$ to $(dd^cu)^{n-1}\wedge T$.
\eT
\par Let $\H_k$ be the Hausdorff measure of dimension $k$. The
second result was proved by Demailly \cite{D2} and then improved by Forn\ae
ss and Sibony \cite[Cor. 3.6]{FS}.
\bT\label{T:fs}
Suppose that $u_0,\dots,u_n$ are plurisubharmonic functions on $D$
and for any choice of indexes $0\le j_1<\dots<j_m\le n$
\begin{equation}\label{e:hmc}
\H_{2(n-m+1)}(L(u_{j_1})\cap\dots L(u_{j_m}))=0.
\end{equation}
Then the function $u_0$ is locally integrable with respect to the
measure $dd^cu_1\wedge\dots dd^cu_n$.
\par Moreover, if for every $0\le j\le n$ the sequences of
plurisubharmonic function $\{u_{jk}\}$ converge to $u_j$ in
$L^1_{loc}(D)$ and $u_{jk}\ge u_j$, then
$$u_{0k}dd^cu_{1k}\wedge\dots\wedge dd^cu_{nk}\to
u_0dd^cu_1\wedge\dots\wedge dd^cu_n$$ in the sense of currents.
\eT
The third theorem we will need is due to Coman \cite[Theorem 3.3]{C}.
\bT\label{T:C}
Let $\Omega$ be a hyperconvex domain in $\aC^n$, $v$ be a
continuous plurisubharmonic exhaustion function on $\Omega$, $T$
be a closed positive (1,1)-current on $\Omega$, and $u$ be a
negative plurisubharmonic function on $\Omega $ such that
$$
\int_\Omega dd^cu\wedge T <\infty.
$$
Then
$$
\int_\Omega vdd^cu\wedge T\geq \int_\Omega udd^cv\wedge T.
$$
\eT

\subsection{Maximal plurisubharmonic functions}

Let $\Omega \subset \aC^n$ be an open set and $u:\Omega \to \rR $
be a plurisubharmonic function. Recall that  $u$ is called {\it
maximal}, if for every relatively compact open subset $G$ in
$\Omega $, and for every upper semicontinuous function $v$ on
$\bar{G}$ such that $v$ is plurisubharmonic in $G$ and $v\leq u$
on $\partial G$, we have $v\leq u$. The following result of
Bedford and Taylor (\cite{Kl},p. 131) gives a necessary and
sufficient condition of maximality.

\bT\label{T:bt} Let $\Omega \subset \aC^n$ be an open set and
$u$ be a plurisubharmonic locally bounded function in $\Omega $.
Then $u$ is maximal if and only if $(dd^cu)^n=0$.
\eT

\par
If $D$ is hyperconvex and $w\in D$, then {\it pluricomplex Green
function with pole at $w$}  is a unique plurisubharmonic in $z$
continuous exhaustion function $g_D(z,w)$ on $D\times D$ such that
$(dd^cg_D(z,w))^n=(2\pi)^n\dl_{w}$ and $|g_D(z,w)-\log|z-w||$ is
bounded on $D$. It is possible to show that
\begin{eqnarray}
g_D(z,w)=\sup \{ u(z): \ u\mbox{ \ is negative and plurisubharmonic in $D$}  \nonumber \\
\mbox{and for some constant $C$,
$u(z)\leq \log |z-w| +C$ near $w$} \}.
\label{eq:Green}
\end{eqnarray}
 (see (\cite{Kl}, p. 221 or \cite{D1} for details).
By Theorem \ref{T:bt}
$g_D(z,w)$ is maximal in $D\setminus \{ w \} $.

\subsection{Lelong-Jensen formula}

\par Let $D$ be a hyperconvex domain in ${\Bbb C}^n$ and $u$
be a continuous negative plurisubharmonic exhausting function
on $D$.
\par We define $B_u(r)=\{z\in
D:\,u(z)< r\}$ and $S_u(r)=\{z\in D:\,u(z)= r\}$. Following \cite{D1} we
let
$$\mu_{u,r}=(dd^cu_r)^n-\chi_{D\sm B_u(r)}(dd^cu)^n,$$
where $u_r=\max\{u,r\}$. The measure $\mu_{u,r}$ is nonnegative and
supported by $S_u(r)$. In \cite[Theorem 1.7]{D1} Demailly had proved the
following fundamental Lelong--Jensen formula.
\bT\label{T:ljf} For all $r<0$ and every plurisubharmonic function
$\phi$ on $D$
$$\mu_{u,r}(\phi)=\il D\phi\mu_{u,r}$$ is finite and
\begin{equation}\label{e:ljf}\begin{align}
\mu_{u,r}(\phi)-&\il {B_u(r)}\phi(dd^cu)^n=\il
{B_u(r)}(r-u)dd^c\phi\wedge(dd^cu)^{n-1}\notag\\&=
\il{-\infty}^rdt\il {B_u(t)}dd^c\phi\wedge(dd^cu)^{n-1}.
\end{align}\end{equation}
\eT
\par The last integral in this formula can be equal to $\infty$.
Then the integral in the left side is equal to $-\infty$. This
cannot happen if $\phi\ge0$.
\par The function
$$\Phi(t)=\il {B_u(t)}dd^c\phi\wedge(dd^cu)^{n-1}<\infty$$ is, evidently,
increasing and by Theorem \ref{T:dem} $\Phi(t)<\infty$ for all
$t<0$. Thus the function
$$\Psi(r)=\il{-\infty}^rdt\il
{B_u(t)}dd^c\phi\wedge(dd^cu)^{n-1}$$ is either identically equal
to $\infty$ or is a continuous function. It follows that  the
function $\mu_{u,r}(\phi)$ is increasing and continuous from the
left.

\subsection{Poincar\'e-Lelong formula}

\par Let $f$ be a holomorphic function on a domain $D\subset \aC^n$
and let $X=\{f=0\}$.
It is well-known that the set $X$ consists of countably many
connected components.  The
critical set $\{\nabla f=0\}$ also consists of countably many
components and the function $f$ is constant on each of components.
Since only finitely many of them intersect a compact set $K\sbs
D$, $f$ has only finitely many critical values on $K$. If $0$ is
in the range of $f$ and is not a
critical value of $f$ ( and, therefore, $f$ is not a constant), then the
analytic set $X$ is a complex manifold of codimension 1. If $0$ is
a critical value of $f$, and $f$ is not a constant function,
 then we denote by $X^{\reg}$ the regular
part of  $X$. The singular part $X^{\sing}$ of $X$ has codimension
at least 2. These and other basic facts about analytic sets
could be found, for example,  in \cite{Chi}.
\par Let $A$ be an irreducible component of $X$. The multiplicity
of $f$ on $A$  is defined as follows (cf. \cite{Chi}). Let $z_0\in
A\cap X^{\reg}$. Since $A$ is a smooth complex manifold, it is
possible to find local coordinates
$\zeta=(\zeta_1,\dots,\zeta_n\}$ near $z_0$ which map $z_0$ to the
point of origin  and the set $A$ into the set $\{\zeta_1=0\}$.
Then in $\zeta $-coordinates we have $f(\zeta)=\zeta_1^mg(\zeta)$,
where $m$ is a natural number and $g(0)\ne0$.
\par Let $h$ be a holomorphic branch of $g^{1/m}$ defined in a
neighborhood of $0$. In coordinates $\xi_1=h(\zeta)\zeta_1$,
$\xi_j=\zeta_j$, $j\ge2$, we have  $f(\xi)=\xi_1^m$.
\par The number $m=m_{A,f}$ is a continuous integer value function
on $X^{\reg}$ and, consequently, is a constant on $A^{\reg}$. It
is called {\it the multiplicity } of $f$ on $A$.
\par  The next theorem contains Poincar\'e--Lelong formula
(\ref{e:plf}). It's proof can be found in \cite[Proposition 3.2.15]{D3}.
\bT\label{T:plf} If $f$ is a holomorphic function on a domain
$D\sbs{\Bbb C}^n$ and $\phi\in\D^{n-1,n-1}(D)$, then
\begin{equation}\label{e:plf}
\il D\phi\wedge dd^c\log|f|=
\sum_Am_{A,f}\il{A^{\reg}}\phi,\end{equation} where the summation
runs over all irreducible components $A$ of $X=\{f=0\}$. Both
sides of this formula define a closed positive current $T$ on $D$
of bidegree $(1,1)$.\eT

\section{Spaces of plurisubharmonic functions}\label{S:spf}

\par Let $D$ be a hyperconvex domain in ${\Bbb C}^n$. We define the space
$PS_u(D)$ as the set of all nonnegative plurisubharmonic functions
$\phi$ on $D$ such that
$$\limsup_{r\to0^-}\mu_{u,r}(\phi)<\infty.$$  (A similar definition
for strongly pseudo-convex domains was given by H{\"o}rmander \cite{H}).
\par Since $\mu_{u,r}(\phi)$ is an increasing function of $r$ for
all $r<0$, we can replace $\limsup$ in the definition of the space
$PS_u(D)$ by $\lim$. So we can introduce the norm on $PS_u(D)$ as
$$\|\phi\|_u=\lim_{r\to0^-}\mu_{u,r}(\phi).$$
\par If ${\Bbb D}$ is the unit disk in ${\Bbb C}$ then the weighted Bergman
space $A^p_\al$ is defined as the set of all holomorphic function
$f$ on ${\Bbb D}$ such that
\begin{equation}\begin{align}
\il{\Bbb D}(1-|z|^2)^\al|f(z)|^p\,dV&=
\il{-\infty}^0\il0^{2\pi}(1-e^{2r})^\al e^r|f(
e^{r+i\th})|^p\,d\th\,dr\notag\\
&\simeq \int_{-\infty}^0 |r|^\alpha  \left (
\int_0^{2\pi}|f(e^{r+i\theta})|^pd\theta \right )dr<\infty.\notag
\end{align}\end{equation}
\par Similar to these spaces we introduce the weighted spaces
$PS_{u,\al}(D)$, $\al>-1$, as the set of all nonnegative
plurisubharmonic functions $\phi$ on $D$ such that
$$\|\phi\|_{u,\al}=\il{-\infty}^0|r|^\al e^r\mu_{u,r}(\phi)\,dr
<\infty.$$ Clearly, $PS_u(D)\sbs PS_{u,\al}(D)\sbs PS_{u,\beta}(D)$ when
$\beta\ge\al$. In what follows, to facilitate the system of notation we let
$PS_u(D)=PS_{u,-1}(D)$ and $\|\phi\|_u=\|\phi\|_{u,-1}$.
\par The following theorem shows that faster decaying near the boundary
of $D$ exhausting functions determine dominating norms.

\bT\label{T:ni} Let $u$ and  $v$ be continuous plurisubharmonic
exhaustion functions on $D$ and let $F$ be a compact set in $D$
such that $F\sbs B_u(r_0)$ for some $r_0<0$ and $v(z)\le u(z)$ for
all $z\in D\sm F$. Then for any $c>1$ and any $a<1-c^{-1}$ we have
$$\mu_{u,r}(\phi)\le c^n\mu_{v,ar}(\phi)$$when
$r\ge r_0$. Moreover, $PS_v(D)\sbs PS_u(D)$ and
$\|\phi\|_u\le\|\phi\|_v$, $PS_{v,\al}(D)\sbs PS_{u,\al}(D)$ and
there is a constant $C$ depending only on $r_0$ and $\al$ such
that $\|\phi\|_{u,\al}\le C\|\phi\|_{v,\al}$.
\eT
\begin{pf} Take a $b>1$ and let $v_1=bv$, $\beta=1-c^{-1}$ and
$\al=c^{-1}$. Take any $r\ge r_0$ and consider the function
$w=\max\{v_1,\al u+\beta r\}$. If $z\in S_u(r)$ then $\al
u(z)+\beta r=r>v_1$ and $w\equiv \al u+\beta r$ on a neighborhood
of $S_{u}(r)$. Moreover, if $z\in S_w(r)$ then $z\in S_u(r)$.
Hence $S_u(r)=S_w(r)$ and $\mu_{w,r}=\al^n\mu_{u,r}$.
\par Let $r_1=\beta r>r$. If $z\in D\sm B_{v_1}(r_1)$, i.e.
$v_1(z)\ge r_1$, then $v_1\ge\beta r>\al u+\beta r$. Hence
$w\equiv v_1$ on a neighborhood of $\ovr{D\sm B_{v_1}(r_1)}$ and
$w$ is an exhaustion function.
\par It follows that
$$b^n\mu_{v,b^{-1}\beta r}(\phi)=\mu_{v_1,\beta r}(\phi)=\mu_{w,r_1}(\phi)
\ge\mu_{w,r}(\phi)=\al^n\mu_{u,r}(\phi).$$ Thus
$$\lim_{b\to1^+}\mu_{v,b^{-1}\beta r}(\phi)\ge\al^n\mu_{u,r}(\phi)$$
and we see that $\mu_{u,r}(\phi)\le c^n\mu_{v,ar}(\phi)$ when $r\ge r_0$.
\par Hence if $\phi\in PS_v(D)$ then $\phi\in PS_u(D)$ and
$\|\phi\|_u\le\|\phi\|_v$.
\par If $\phi\in PS_{v,\al}(D)$ then we let $c=2$ and $a=1/4$. We write
$$\il{-\infty}^0|r|^\al e^r\mu_{u,r}(\phi)\,dr=
\il{-\infty}^{r_0}|r|^\al e^r\mu_{u,r}(\phi)\,dr+\il{r_0}^0|r|^\al
e^r\mu_{u,r}(\phi)\,dr.$$ The first integral does not exceed
$C_1(r_0,\al)\mu_{u,r_0}(\phi)$, where
$$C_1(r_0,\al)=\il{-\infty}^{r_0}|r|^\al e^r\,dr.$$
As we proved above $\mu_{u,r_0}(\phi)\le2^n\mu_{v,r_0/4}(\phi)$, and
$$\|\phi\|_{v,\al}\ge \il{r_0/4}^0|r|^\al
e^r\mu_{v,r}(\phi)\,dr\ge \mu_{v,r_0/4}(\phi)\il{r_0/4}^0|r|^\al e^r\,dr.$$
Hence,
$$\il{-\infty}^{r_0}|r|^\al e^r\mu_{u,r}(\phi)\,dr\le
\frac{2^nC_1(r_0,\al)}{C_2(r_0,\al,c)}\|\phi\|_{v,\al},$$ where
$$C_2(r_0,\al,c)=\il{r_0/4}^0|r|^\al e^r\,dr.$$
\par The second integral does not exceed
$$2^n\il{r_0}^0|r|^\al e^r\mu_{v,r/4}(\phi)\,dr=
2^{2\al+n+2}\il{r_0/4}^0|t|^\al e^{4t}\mu_{v,t}(\phi)\,dt,$$ which, in turn,
does not exceed
$$2^{2\al+n+2}\il{r_0/4}^0|t|^\al
e^t\mu_{v,t}(\phi)\,dt\le 2^{2\al+n+2}\|\phi\|_{v,\al}.$$ Combining these
estimates we get that
$$\|\phi\|_{u,\al}\le\left(\frac{2^nC_1(r_0,\al)}{C_2(r_0,\al,c)}+
2^{2\al+n+2}\right)\|\phi\|_{v,\al}.$$ We are done.
\end{pf}
\par This theorem has a couple of useful corollaries.
\bC\label{C:ni}Let $u$ and  $v$ be continuous plurisubharmonic
exhaustion functions on $D$ and let $F$ be a compact set in $D$
such that $bv(z)\le u(z)$ for some constant $b>0$ and all $z\in
D\sm F$. Then $PS_v(D)\sbs PS_u(D)$ and $\|\phi\|_u\le
b^n\|\phi\|_v$, $PS_{v,\al}(D)\sbs PS_{u,\al}(D)$ and there is a
constant $C$ depending only on $r_0$, $b$ and $\al$ such that
$\|\phi\|_{u,\al}\le C\|\phi\|_{v,\al}$.
\eC
\begin{pf}  Remark that $bv_r=\max\{bv,br\}$, and, thus,
$b^n\mu_{v,r}=\mu_{bv,br}$. Now  Corollary follows immediately from Theorem
\ref{T:ni}
\end{pf}
\bC\label{C:es}Let $u$ and $v$ be continuous plurisubharmonic
exhaustion functions on $D$ and let $F$ be a compact set in $D$
such that
\begin{equation}\label{e:ef}
bv\le u\le b^{-1}v
\end{equation}
for some constant $b>0$ and all $z\in D\sm F$. Then $PS_v(D)=
PS_u(D)$, $PS_{v,\al}(D)=PS_{u,\al}(D)$ and the identity mappings
are continuous.
\eC
\par Let us describe a class of functions for which the inequality
(\ref{e:ef}) holds automatically. We denote  by $\E_0(D)$ the set of all
plurisubharmonic exhaustion functions $u$ on $D$ with compactly supported
$(dd^cu)^n$. The set $\E_0(D)$ is not empty: pluricomplex Green functions
are in $\E_0(D)$. Moreover, by (\ref{e:ljf}) if $u\in\E_0(d)$ then the
space $A^p_{u,\al}(D)$ contains constants and, consequently, all bounded
holomorphic functions.
\bL\label{L:cd} If $u,v\in\E_0(D)$ then there is a constant $b>0$
such that $bu(z)\le v(z)\le b^{-1}u(z)$ near $\bd D$.\eL
\begin{pf} Let us take a compact set $F\sbs D$
containing the supports of $(dd^cu)^n$ and $(dd^cv)^n$ such that both $u$
and $v$ are bounded on $\bd F$ and (\ref{e:ef}) holds on $\bd F$ for some
number $b>0$. By the maximality of $u$ and $v$ on $D\sm F$ this inequality
holds on $D\sm F$ also.
\end{pf}
\par Thus functions in $\E_0(D)$ generate the same spaces $PS_{u,\al}(D)$
and $PS_u(D)$ with equivalent norms. As the following proposition shows
these spaces are the largest in our class.
\bP\label{P:ls}Let $u\in\E_0(D)$ and let $v$ be a continuous plurisubharmonic exhaustion
function on $D$. Then $PS_{v,\al}(D)\sbs PS_{u,\al}(D)$ for each $\al\ge-1$
and there is a constant $C(\al)$ such that $\|\phi\|_{u,\al}\le
C(\al)\|\phi\|_{v,\al}$.
\eP
\begin{pf} Take a number $r<0$ such that $\supp(dd^cu)^n\sbs B_u(r)$.
There is a constant $b>0$ such that $bv\le u$ on $S_u(r)$. By the
maximality of $u$ the same inequality holds on $D\sm B_u(r)$. Now the
proposition follows from Corollary \ref{C:ni}.
\end{pf}
\par Using pluricomplex Green functions we can get estimates for
point evaluations. If $u(z)=g_D(z,w)$, $w\in D$, then by (\ref{e:ljf}) (see
also \cite[Thm. 5.1]{D1}) we get
\begin{equation}\label{e:pe}(2\pi)^n\phi(w)=
\mu_{u,r}(\phi)-\il
{B_u(r)}(r-u)dd^c\phi\wedge(dd^cu)^{n-1}\le\mu_{u,r}(\phi).\end{equation}

\bT\label{T:ci} Let $v$ be a continuous plurisubharmonic
exhaustion function on $D$. Then for any compact set $K\sbs D$ and any
$\al\ge-1$ there is a constant $c$ such that for all $w\in K$ and all
nonnegative plurisubharmonic functions on $D$ we have $\phi(w)\le
c\|\phi\|_{v,\al}$.
\eT

\begin{pf} By \cite[Th. 4.14]{D1} the function $u_w(z)=g_D(z,w)$ is
continuous on $D\times D$. If $V\sbs\sbs D$ is an open set
containing $K$, then there is a negative constant $a$ such that
$a\le u_w(z)$ for all $w\in K$ and $z\in\bd \ovr V$. Consequently,
there is a constant $c>0$ such that $v(z)\le cu_w(z)$ for all
$w\in K$ and $z\in\bd V$. The maximality of $g_D$ outside $F$
implies that that $v(z)\le cu_w(z)$ for all $w\in K$ and $z\in
D\sm\ovr V$.
\par There is a number $r_0<0$ such that $\ovr V\sbs
B_{u_w}(r_0)$ for all $w\in K$. Hence  by (\ref{e:pe}) and Theorem
\ref{T:ni}
$$\phi(w)\le\frac{c^n}{(2\pi)^n}\mu_{u_w,r}(\phi)\le
\frac{c^n}{(2\pi)^n}\|\phi\|_{v,-1}$$ and
$$\phi(w)\il{r_0}^0|r|^\al e^r\,dr\le
\frac{c^n}{(2\pi )^n}\il{-\infty}^0|r|^\al e^r\mu_{u_w,r}(\phi)\,dr\le
C\|\phi\|_{v,\al}.$$
\end{pf}
\par We denote by $\E(D)$ the class of continuous
plurisubharmonic exhaustion functions $v$ on $D$ such that the inequality
$bv\le u\le b^{-1}v$ holds near $\bd D$ for some function $u\in\E_0(D)$ and
constant $b>0$. By Corollary \ref{C:es} these exhaustion functions generate
the same space as functions from $\E_0(D)$ with equivalent norms. The class
$\E(D)$ contains many important functions.
\bP\label{P:sc} A plurisubharmonic exhaustion function $v\in C^1(\ovr D)$ on a
domain $D$ belongs to $\E(D)$.\eP
\begin{pf} Fix $w\in D$ and let $u(z)=g_D(z,w)$. As it was
explained in the proof of Proposition \ref{P:ls} the inequality $bv\le u$
holds for some constant $b$ near $\bd D$. On the other hand by Hopf's lemma
$u(z)\le-cd(z,\bd D)$, where $c$ is some positive constant and $d(z,\bd D)$
is the distance from $z$ to the boundary of $D$. Since for some positive
constant $a$ the function $\rho(z)>-ad(z,\bd D)$ near $\bd D$, the
proposition follows.
\end{pf}
\par The classical approach to the definition of spaces
$PS_{u,\al}(D)$ and $PS_u(D)$  would be to restrict the studies to
domains $D$ for which there is a smooth function $\rho$ defined on
a neighborhood of $\ovr D$ and such that $D=\{\rho<0\}$ and
$\nabla\rho\ne0$ on $\bd D$. Then the definitions of spaces will
go through as above with the replacement of $\mu_{u,r}(\phi)$ by
$$\nu_{\rho,r}(\phi)=\il{\{\rho=r\}}\phi\,d\si,$$
where $d\si$ is the surface measure. We will denote the new norms
by $\|\phi\|_\rho$ and $\|\phi\|_{\rho,\al}$.
\par The following theorem shows that both approaches coincide in
the case of strongly pseudoconvex domains.
\bT\label{T:cl} Let $D$ be a strongly pseudoconvex domain and let
$\rho$ be a strictly plurisubharmonic function defined on a neighborhood of
$\bar D$ such that $D=\{\rho<0\}$ and $\nabla\rho\ne0$ on $\bd D$. If
$u\in\E(D)$ then for each $\al\ge-1$ there is a constant $C>1$ such that
$C^{-1}\|\phi\|_{u,\al}\le\|\phi\|_{\rho,\al}\le C\|\phi\|_{u,\al}$.
\eT
\begin{pf} We just note that by \cite[(1.5)]{D1}
$\mu_{\rho,r}=(dd^c\rho)^{n-1}\wedge d^c\rho$. Since this form is
continuous and strictly positive on $D$, there are positive constants $a_1$
and $a_2$ such that $a_1\mu_{\rho,r}\le d\si\le a_2\mu_{\rho,r}$. The rest
follows from Proposition \ref{P:sc} and Corollary \ref{C:es}.
\end{pf}

\section{Hardy and Bergman spaces}\label{S:HBs}

\par Let $u$ be a continuous plurisubharmonic exhaustion
function on a hyperconvex domain $D$. We define the Hardy space $H^p_u(D)$,
$p>0$, as the space of all holomorphic functions $f$ on $D$ such for
$|f|^p\in PS_u(D)$. The Hardy norm $\|f\|_{H^p_u}=\||f|^p\|^{1/p}_u$.
\par We define the Bergman space $A^p_{u,\al}(D)$, $p>0$, as the
space of all holomorphic functions $f$ on $D$ such for $|f|^p\in
PS_{u,\al}(D)$. The Bergman norm
$\|f\|_{A^p_{u,\al}}=\||f|^p\|^{1/p}_{u,\al}$.
\par Similar to the case of nonnegative plurisubharmonic functions we
denote $A^p_{u,-1}(D)=H^p_u(D)$ and $\|f\|_{A^p_{u,-1}}=\||f|^p\|_u$.
\par The classical definition of Hardy and Bergman spaces when
$D$ is the unit ball $B\sbs{\Bbb C}^n$ instead of the measures
$\mu_{u,t}$ uses the measures $\nu_r=dS$, where $dS$ is the
normalized surface area, on spheres of radius $r$. If we take
$u(z)=\log|z|$ as an exhaustion function, then $dd^cu_t$ is a
rotationally invariant measure which is supported by the sphere
$\{|z|=r=e^t\}$ and is a constant multiple of $dS$. By Riesz
representation formula or by the Lelong--Jensen formula
$$\il {\Bbb D}dd^cu_t=(2\pi)^n.$$ So $\mu_{u,t}=(2\pi)^ndS$ and
we see that our definition coincides with the classical one.
\par If $D={\Bbb D}^2$ is the unit bidisk in ${\Bbb C}^2$ with
coordinates $z_1$ and $z_2$ and $u(z)=\log\max\{|z_1|,|z_2|\}$,
then $(dd^cu)^2=(2\pi)^2\dl_0$ and the measure $(dd^cu_r)^2$ is
supported by the set $\{|z_1|=|z_2|=r\}$. The latter measure is
also rotationally invariant, so $\mu_{u,r}$ is a constant multiple
of the normalized surface area on the torus $\{|z_1|=|z_2|=r\}$
and again our definition coincides with the classical definition
of Hardy (but not Bergman) spaces.
\par Let us show that Hardy and Bergman spaces are Banach.
\bT\label{T:hs} Let $u$ be a continuous plurisubharmonic
exhaustion function on $D$. Then the spaces $A^p_{u,\al}(D)$, $p\ge1$, are
Banach.
\eT
\begin{pf} Clearly $\|u\|_{A^p_{u,\al}}$ is a norm on $A^p_{u,\al}(D)$.
\par Suppose that $\{f_j\}$ is a Cauchy sequence in $H^p_u(D)$.
By Theorem \ref{T:ci} $\{f_j\}$ is a Cauchy sequence in the
uniform metric on any compact set in $D$. Hence this sequence
converges to a holomorphic function $f$ on $D$ uniformly on
compacta. In particular, for every fixed $r<0$
$\mu_{u,r}(|f-f_j|^p)\to 0$ as $j\to\infty$.
\par There is $A>0$ such that $\|f_j\|_{H^p_u}\le A$ for all $j$. Hence
$$\left(\mu_{u,r}(|f|^p)\right)^{1/p}\le
\left(\mu_{u,r}(|f-f_j|^p)\right)^{1/p}+A$$ and we see that $f\in
H^p_u(D)$ and $\|f\|_{H^p_u}\le A$.
\par Suppose that $\limsup\|f-f_j\|_{H^p_u}\ne0$. Switching to a
subsequence we may assume that $\|f-f_j\|_{H^p_u}\ge a>0$ for all
$j$. Take $i$ such that $\|f_i-f_j\|_{H^p_u}<a/4$ when $j>i$
and then find $r<0$ such that $\mu_{u,r}(|f-f_i|^p)>(a/2)^p$. Then
$$\left(\mu_{u,r}(|f-f_j|^p)\right)^{1/p}\ge
\left(\mu_{u,r}(|f-f_i|^p)\right)^{1/p}-\|f_i-f_j\|_{H^p_u}\ge\frac
a4.$$ But the left side converges to 0 as $j\to\infty$ and this
proves that $\lim\|f-f_j\|_{H^p_u}=0$.
\par The proof for the Bergman spaces is basically the same.
\end{pf}
\par Our next goal is to prove a multidimensional analog of
Littlewood's subordination  principle.
\bT\label{T:lmp} Let $D$ be a hyperconvex domain, $z_0\in D$,
$u(z)=g_D(z,z_0)$ and let $f:\,D\to{\Bbb D}$ be a holomorphic
mapping with $f(z_0)=w_0$. If $v(w)=\log|(w-w_0)/(1-\ovr w_0w)|$
then $\mu_{u,r}(\phi\circ f)\le(2\pi)^{n-1}\mu_{v,r}(\phi)$ for
every subharmonic function $\phi$ on ${\Bbb D}$.
\eT
\begin{pf} The relation (\ref{eq:Green}) implies that  $v^*(z)=v(f(z))\le u(z)$ on $D$. Hence
$f(B_u(r))\sbs B_v(r)$.
\par Let $r<t<0$ and let $h_t$ be the harmonic function on
$B_v(t)$ equal to $\phi$ on $S_v(t)$. Then
$h_t(f(z))\ge\phi(f(z))$ on $S_u(r)$ and therefore
$$\mu_{u,r}(\phi\circ f)\le\mu_{u,r}(h_t\circ f).$$
Since the function $h_t\circ f$ is pluriharmonic, $dd^c(h_t\circ
f)\equiv0$. Since $h_t\circ f$ is defined on $B_u(t)$ we can use
(\ref{e:pe}) to get
$$\mu_{u,r}(h_t\circ f)=(2\pi)^nh_t(f(z_0))=(2\pi)^nh_t(w_0).$$
Since $\mu_{v,t}(\phi)=\mu_{v,t}(h_t)=2\pi h_t(w_0)$, we see that
$\mu_{u,r}(\phi\circ f)\le\mu_{v,t}(\phi)$. But
$\mu_{v,r}(\phi)=\lim_{t\to r^+}\mu_{v,t}(\phi)$ and the the
result follows.
\end{pf}
\par The following corollary follows immediately from the
previous theorem.
\bC\label{C:1d} Let $v(w)=\log|w|, \ w\in {\Bbb D}$, and
let $D$ be a hyperconvex domain, $u\in\E(D)$,  and  $f:\,D\to{\Bbb D}$
be a holomorphic mapping. Then the composition operator $C_f$ maps
$A_{v,\al}({\Bbb D})$ into $A_{u,\al}(D)$. Moreover, there is a constant
$A$ depending only on $|w_0|$, $w_0=f(z_0)$, and $u$ such that
$\|C_f\phi\|_{u,\al}\le A\|\phi\|_{v,\al}$.
\par Consequently, $C_f$ maps continuously $A^p_{v,\al}({\Bbb D})$ into $A^p_{u,\al}(D)$.
\eC
\par For a holomorphic mapping of a domain $D$ into another domain $F$
we introduce the measure $\nu_{f,\al}$ on $F$ defined as
$$\nu_{f,\al}(E)=
\il{-\infty}^0|r|^\al e^r\mu_{u,r}(\chi_E(f(z)))\,dr ,$$ where
$\chi_E$ is the characteristic function of $E$. The proof of the
following lemma is similar to  the proof of Theorem 2.6 in
\cite{CoM}.
\bL\label{L:n1cm} Let $E(w,s)=\{|z-w|<s\}\cap{\Bbb D}$.
In the assumptions of Corollary \ref{C:1d} there is a constant
$a\ge1$ depending only on $|w_0|$ such that
$\nu_{f,\al}(E(w,s))\le as^{\al+2}$ when $s\le1$ for all $w\in S$.
\eL
\begin{pf} It suffices to prove this lemma for $w=1$. Let us take
$t=1/(1+2s)$, $0<s<1$, and consider the function
$$\phi(z)=\frac1{|1-tz|^{2\al+4}}.$$
Then by Corollary \ref{C:1d}
\begin{equation}\begin{align}&\il{\Bbb
D}\phi\,d\nu_{f,\al}=\il{-\infty}^0|r|^\al
e^r\mu_{u,r}(\phi(f(z)))\,dr \notag\\&=\|\phi\circ
f\|_{A^1_{u,\al}(D)}\le c\|\phi\|_{A^1_{v,\al}({\Bbb
D})}=c_1\frac1{(1-t^2)^{\al+2}}.\end{align}\end{equation}
\par If $|1-z|<s$, then
$$\frac1{|1-tz|}=\frac{1+2s}{|1-z+2s|}\ge\frac{1+2s}{3s}\ge
\frac1{3s}$$ and
$$\frac1{1-t^2}\le\frac3s.$$
Thus
$$c_1\left(\frac3s\right)^{\al+2}\ge c_c\frac1{(1-t^2)^{\al+2}}\ge
\il{\Bbb D}\phi\,d\nu_{f,\al}\ge
\left(\frac1{3s}\right)^{2\al+4}\nu_{f,\al}(E(w,s)).$$ Hence
$\nu_{f,\al}(E(w,s))\le as^{\al+2}$ for all $w\in S$.
\end{pf}

\section{The change of variables formula}\label{S:cvf}

\par Throughout this section $D$ is a domain in ${\Bbb C}^n$ and
$u_0,\dots,u_{n-1}$ are continuous plurisubharmonic functions on
$D$. We set $$T=u_0dd^cu_1\wedge\dots\wedge dd^cu_{n-1}.$$ We will
assume that for any choice of indexes $0\le j_1<\dots<j_m\le n-1$
\begin{equation}\label{e:hmct}
\H_{2(n-m+1)}(L(u_{j_1})\cap\dots \cap L(u_{j_m}))=0.
\end{equation}
We fix a holomorphic function $f$ on $D$ mapping $D$ into a domain
$\Om\sbs{\Bbb C}$. For a subharmonic function $v$ on $\Om$ we define
 a functional on $\D^{0,0}(D)\times\D^{0,0}(\Om)$
$$\beta_{f,v,T}(\phi,\psi)=\beta_{f}(\phi,\psi)=
\il D\phi\psi^*dd^cv^*\wedge T,$$ where $\psi^*$ and $v^*$ are
$\psi\circ f$ and $v\circ f$ respectively.
\bL\label{L:ncf} If the functions $u_0,\dots,u_{n-1}$ satisfy
(\ref{e:hmct}) and for any choice of indexes $0\le
j_1<\dots<j_m\le n-1$
$$\H_{2(n-m)}(L(v^*)\cap L(u_{j_1})\cap\dots \cap L(u_{j_m}))=0,$$ then
for a fixed $\phi$ the functional $\beta_{f,v,T}(\phi,\psi)$ is a
current of bidegree $(0,0)$ on $\Om$ and for a fixed $\psi$ the
functional $\beta_{f,v,T}(\phi,\psi)$ is a current of bidegree
$(0,0)$ on $D$.
\eL
\begin{pf} It follows from Theorem \ref{T:fs} that the current
$dd^cv^*\wedge T$ has a locally finite mass. Since the function
$\psi^*$ is bounded on $D$, we see that $\beta_{f,v,T}(\phi,\psi)$
is a current of bidegree $(0,0)$ on $D$.
\par Let $\phi\in\D^{0,0}(D)$ be fixed. Then the functional
$\beta_{f,v,T}(\phi,\psi)$ is defined for all
$\psi\in\D^{0,0}(\Om)$. To show that the functional
$\beta_{f,v,T}(\phi,\psi)$ is continuous we note that if a
sequence $\{\psi_j\}$ converges to 0, then the functions
$\{\phi\psi_j\}$ converge to 0 uniformly on $D$ and, therefore,
$\beta_{f,v,T}(\phi,\psi_j)\to0$ as $j\to\infty$.
\end{pf}
\par The main goal of this section is to get a formula for the
measure on $\Om$ defining $\beta_{f,v,T}$. For this we take a
function $\phi\in\D^{0,0}(D)$ and for $w\in{\Bbb C}$ introduce the
function
$$N_{f,T}(w,\phi)=N(w,\phi)=\il D\phi dd^c\log|f-w|\wedge T.$$
\par This function need not  be finite. For example if
$D={\Bbb D}$, $T=\log|z|$, $f(z)=z$ and $w=0$, then
$N(0,\phi)=-\infty$ when $\phi(0)>0$. The following proposition
gives sufficient conditions for the function $N(w,\phi)$ to be
finite and states some of its properties.
\par  Let $X_w=\{f=w\}$. Write
$$N_0(w,\phi)=
\sum_Am_{A,f}\il{A^{\reg}}\phi T,$$ where the summation runs over all
irreducible components $A$ of $X_w$ and $m_{A,f}$ is the multiplicity of
$f$ on $A^{\reg}$.
\bP\label{P:ln} If the
functions $u_0,\dots,u_{n-1}$ satisfy (\ref{e:hmct}) and
$L(u_j)\cap X_w=\emptyset$ for any $0\le j\le n-1$, then
$N(w,\phi)<\infty$ and $N(w,\phi)=N_0(w,\phi)$. Moreover, if,
additionally, $\phi\ge0$ on $D$ and the function $u_0<0$ on
$\supp\phi$, then the function $N(\tau,\phi)$ is upper
semicontinuous (as a function of $\tau$)at $w$. Moreover, if the
set $X_w$ is smooth then $N(\tau,\phi)$ is continuous at $w$. \eP
\begin{pf} The fact that $N(w,\phi)$ is bounded follows immediately
from Theorem \ref{T:fs}.
\par To show that $N(w,\phi)=N_0(w,\phi)$ we observe that the function
$u_0$ does not exceed some number $a>0$ on a neighborhood
$V\sbs\sbs D$ of $\supp\phi$. Thus
$$T=(u_0-a)dd^cu_1\wedge\dots\wedge dd^cu_{n-1}+
add^cu_1\wedge\dots\wedge dd^cu_{n-1}$$ is the sum of a negative
and a positive current. Hence, without loss of generality we may
assume that $u_0<0$ on $V$.
\par Then we take a decreasing sequence of smooth
plurisubharmonic functions $u_{jk}$ defined on $V$ and converging
to $u_j$ for every $0\le j\le n-1$. By Theorem
\ref{T:fs}
$$\il {D}\phi T\wedge
dd^c\log|f-w|=\lim_{k\to\infty}\il {V}\phi T_k\wedge
dd^c\log|f-w|,$$ where $T_k=u_{0k}dd^cu_{1k}\wedge\dots\wedge
dd^cu_{n-1,k}$.
\par By Poincar\'e--Lelong formula (\ref{e:plf})
$$\il {V}\phi T_k\wedge dd^c\log|f-w|=
\sum_Am_{A,f}\il{A^{\reg}}\phi T_k.$$
\par Now we choose an increasing sequence of nonnegative functions
$\psi_j\in\D^{0,0}(V)$ equal to 0 on $X^{\sing}_w$ and converging
to 1 on $\supp\phi\sm X^{\sing}_w$. Since the functions $u_{jk}$
are bounded on $X_w$ and converge to $u_j$,
$$\lim\limits_{k\to\infty}\il{A^{\reg}}\phi\psi_jT_k=
\il{A^{\reg}}\phi\psi_jT.$$ Thus,
$$\il {V}\phi\psi_jT\wedge
dd^c\log|f-w|= \sum_Am_{A,f}\il{A^{\reg}}\phi\psi_jT.$$
\par If $\phi\ge0$ on $V$, then letting $j\to\infty$ by the
monotone convergence theorem we get that
\begin{equation}\label{e:le}
\il {V\sm X^{\sing}_w}\phi T\wedge dd^c\log|f-w|=
\sum_Am_{A,f}\il{A^{\reg}}\phi T.
\end{equation}
If $\phi$ is not positive on $V$ then we take a function
$\wtl\phi\in\D^{0,0}(V)$ such that $\wtl\phi\ge|\phi|$. Since the
function $\psi=\wtl\phi-\phi\ge0$ on $D$ and $\phi=\wtl\phi-\psi$,
the equality (\ref{e:le}) holds in this case also.
\par Now let us show that
$$\il {X^{\sing}_w}\phi T\wedge dd^c\log|f-w|=0.$$
Using a partition of unity we can reduce it to the case when $D=
B(z_0,r)$, $\phi\ge0$, $z_0\in X^{\sing}_w$ and the analytic set
$X^{\sing}_w$ is the set of common zeros of holomorphic function
$f_1,\dots,f_k$ on $D$. Then for the plurisubharmonic function
$$\wtl u_0=u_0+\log\sum_{i=1}^k|f_i|^2$$ the set
$L(\wtl u_0)=L(u_0)\cup X^{\sing}_w$. If $\wtl T=\wtl
u_0dd^cu_1\wedge\dots\wedge dd^cu_{n-1}$ and $u_n=\log|f-w|$, then the
conditions (\ref{e:hmct}) hold when al least one of the indexes $j_i$ is
between 1 and $n-1$ because $L(u_j)\cap X_w=\emptyset$ for all $j$. If
$m=1$ then $\H_{2n}(L(u_0))=\H_{2n}(X_w)=0$, and, therefore,
$\H_{2n}(L(\wtl u_0 ))=0$.
 If $m=2$, $j_1=0$ and $j_2=n$,
then the set $L(\wtl u_0)\cap L(u_n)=X^{\sing}_w$ is an analytic
set of codimension at least 2 and, therefore, $\H_{2n-2}(L(\wtl
u_0)\cap L(u_n))=0$. Hence, by Theorem \ref{T:fs}, the current
$\wtl T\wedge\log|f-w|$ is locally integrable. But $\wtl
u_0\equiv-\infty$ on $X^{\sing}_w$. Therefore, the current
$dd^cu_1\wedge\dots dd^cu_n$ has no mass on $X^{\sing}_w$, and
$$\il {X^{\sing}_w}\phi T\wedge dd^c\log|f-w|=0.$$
\par Let us show that the function $N_0(w,\phi)$ is upper
semicontinuous at $w$ which we assume to be equal to 0. Take a
compact set $K_A$ in a component $A^{\reg}_0$, cover $K_A$ by open
balls $V_k\sbs V$, $1\le k\le k_0$, such that
$f(\zeta)=\zeta_1^m$, $m=m_{A,f}$, in some new coordinates
$(\zeta_1,\dots,\zeta_n\}$ on $V_k$ and  form a partition of unity
$\psi_k\in\D^{0,0}(V_k)$ such that $\psi=\sum_j\psi_k=1$ on a
neighborhood $W_A$ of $K_A$.
\par There is an $\eps>0$ such that for every $k$ and every $w\ne0$,
$|w|<\eps$, the set $X_w\cap V_k$ is a complex manifold consisting
of $m$ connected components $X^w_{kj}$, $1\le j\le m$, given by
the equations $\zeta_1=w^{1/m}$.
\par If $'\zeta=(\zeta_2,\dots,\zeta_m)$ then by continuity of the
functions $u_i$ the functions
$u^w_{ikj}('\zeta)=u_i(w^{1/m},'\zeta)$ uniformly converge to
$u_i(0,'\zeta)$ as $w\to0$ and, therefore,
$$\lim\limits_{w\to0}\il{X^w_{kj}}\phi\psi_kT=
\il{X_0\cap V_k}\phi\psi_kT.$$ Since $\phi\ge0$ and $u_0\le0$, so
$\phi T$ is a negative current, summation over $j$ and $k$ results
in
$$\il{K_A}\phi T\ge \limsup\limits_{w\to0}\il{X_w\cap W_A}\phi T,$$
and our statement follows.
\par To finish the proof of the proposition let us show that
$N_0(0,r)$ is continuous at $0$ when $X_{0}$ is smooth. In this
case we can take a compact set $K=X_0\cap\supp\phi$ and the balls
$V_k$, covering $K$, can be chosen so that each of them contains
only one component of $X_0$. Now the previous argument yields the
proof.
\end{pf}
{\bf Remark.} This result does not hold when $L(u_j)\cap X_w\ne\emptyset$.
For example, if $D$ is the unit ball in ${\Bbb C}^2$ with coordinates $z_1$
and $z_2$, $u_0(z)\equiv1$, $u_1(z)=\log|z|$ and $f(z)=z_1z_2$, then
$$N(0,\phi)=\il {D}\phi dd^cu_1\wedge dd^c\log|f|=2\phi(0).$$
However, $u_1(z)$ is equal to either $\log|z_1|$ or $\log|z_2|$ on
$X^{\reg}_0$. Therefore,  $dd^cu_1=0$ on $X^{\reg}_0$ and
$N_0(0,\phi)=0$.
\par We will need the following version of Fubini's theorem.
\bL\label{L:dp} Let $X$ be a complex manifold of dimension $n-1$,
$U$ be a domain in ${\Bbb C}$, $v$ be a subharmonic function on
$U$ and  $u_0,\dots,u_{n-1}$ be bounded plurisubharmonic
functions on $X\times U$. If $\phi\in\D^{0,0}(X\times U)$ then
$$\il{X\times U}\phi T\wedge dd^cv=
\il{U}\left(\il{X\times\{w\}}\phi T\right)dd^cv.$$
\eL
\begin{pf} Let us show that the lemma holds when the functions $u_j$
and $v$ are smooth. Introducing a partition of unity we can reduce
it to the case when $X$ is a domain in ${\Bbb C}^{n-1}$ with
coordinates $z=(z_1,\dots,z_{n-1})$. Let $w$ be a coordinate on
$U$ and
$$dd_z^cu_j=
2i\sum_{i,k=1}^{n-1}\frac{\bd^2 u_j}{\bd z_i\bd\ovr
z_k}(z,w)dz_i\wedge d\ovr z_k.$$ Note that $dd^c_zu_j$ is the
restriction of $dd^cu_j$ to $X\times\{w\}$. Hence
$$\il{X\times\{w\}}\phi
T=\il{X\times\{w\}}\phi T_w,$$ where
$T_w=u_0dd^c_zu_1\wedge\dots\wedge dd^c_zu_{n-1}$. Since $\phi
T\wedge dd^cv=\phi T_w\wedge dd^cv$, the lemma follows immediately
from Fubini's theorem.
\par In the general case let us fix  open sets
$U_1\sbs\sbs U$ and $X_1\sbs\sbs X$ such that $\supp\phi\sbs
X_1\times U_1$. We take decreasing sequences $\{u_{jk}\}$, $0\le
j\le n-1$, and $\{v_k\}$ of smooth plurisubharmonic functions on
$X_1\times U_1$ and $U_1$ converging to $u_j$ and $v$
respectively.
\par Let $T_k=u_{0k}dd^cu_{1k}\wedge\dots\wedge u_{n-1,k}$. Note
that the functions
$$\Phi_k(w)=\il{X\times\{w\}}\phi T_k$$ are
smooth and have compact support on $U$. Since $dd^cv_m$ converge
weak-$*$ to $dd^cv$ on $U$ we see that
$$\il{X\times U}\phi T_k\wedge dd^cv=
\il{U}\Phi_{k}(w)\,dd^cv.$$
\par As $k\to\infty$ the monotonic convergence of $u_{jk}$ implies
that
$$\lim_{k\to\infty}
\il{X\times U}\phi T_k\wedge dd^cv=\il{X\times U}\phi T\wedge
dd^cv$$ and for each $w$ the functions $\Phi_k(w)$ converge to
$$\Phi(w)=\il{X\times\{w\}}\phi T.$$
\par Since $u_{jk}\ge u_{j,k+1}\ge u_j$ and the
functions $u_j$ are bounded, there is a constant $a$ such that the
$L^\infty$-norms of $u_{jk}$ on $X_1\times U_1$ do not exceed $a$.
By Chern-Levine-Nirenberg inequality (see \cite{Kl}, p.111-112) there is
another constant $C$ such that $|\Phi_k(w)|\le C$. The support of
$\Phi_k(w)$ lies in $U_1\sbs\sbs U$ and $dd^cv$ has a finite mass
on $U_1$. Thus, by the dominated convergence theorem
$$\lim_{k\to\infty}\il {U}\Phi_k(w)dd^cv=\il {U}
\Phi(w)dd^cv.$$ Hence,
$$\il{X\times U}\phi T\wedge dd^cv=
\il {U}\Phi(w)dd^cv.$$
\end{pf}
\par The following theorem, which is the main goal of this section,
gives us a change of variables formula.

\bT\label{T:cvi}Let $f$ be a holomorphic function on a domain
$D\sbs{\Bbb C}^n$ and  $u_0,\dots,u_{n-1}$ be continuous
plurisubharmonic functions on $D$. Let $\Om$ be a domain in ${\Bbb
C}$ such that $f(D)\sbs\Om$ and let $v$ be a subharmonic function
on $\Om$ and $v^*=v\circ f$. Suppose that the sets $L(u_j)$ are
finite and $f(L(u_j))\cap L(v)=\emptyset$ for $0\le j\le n-1$. If
$\phi\in\D^{0,0}(D)$ and $\psi\in\D^{0,0}(\Om)$, then
$$\beta_{f,v,T}(\phi,\psi)=
\il{\Om} N(w,\phi\psi^*)dd^cv.$$
\eT

\begin{pf} Let us assume for a while that the functions $\phi$ and
$\psi$ are nonnegative and $u_0\le0$ on $\supp\phi$. For a set
$E\sbs D$ we define the measure
$$\nu(E)=
-\il{E}\phi\psi^* T\wedge dd^cv^*.$$
\par Let $W_1$ be the union of the sets $f(L(u_j))$, $0\le j\le n-1$.
This is a finite set consisting of $p$ points $\{w_j\}$. Let
$X=\cup_{j=1}^pX_{w_j}$ and $\wtl
u_0=u_0+\sum_{j=1}^p\log|f-w_j|$. Let $u_n=v^*$. If $j_m=n$ and
$0\le j_k\le n-1$ at least once for a choice of indexes $0\le
j_1<\dots<j_m$, then $L(u_{j_1})\cap\dots\cap
L(u_{j_m})=\emptyset$ and
$$\H_0(L(u_{j_1})\cap\dots\cap L(u_{j_m}))=0.$$  If $m\ge 2$ and
$1\le j_m<n$ then the sets $L(u_{j_1})\cap\dots\cap L(u_{j_m})$
are finite and $\H_1(L(u_{j_1})\cap\dots\cap L(u_{j_m}))=0$. If
$m=1$ and $j_1=0$ or $j_1=n$ then $\H_{2n}(X)=0$. So we see that
(\ref{e:hmc}) holds and the current $\wtl
u_0dd^cu_1\wedge\dots\wedge dd^cu_{n-1}\wedge dd^cv^*$ has locally
finite mass. Since $\wtl u_0=-\infty$ in $X$,
$$\il {X}dd^cu_1\wedge\dots\wedge dd^cu_{n-1}\wedge dd^cv^*=0$$ and,
consequently, $\nu(X)=0$.
\par Now, let $W_2=\{w_1,\dots,w_q\}$ be the set of those critical
values of $f$ on $\supp\phi$ which do not belong to any
$f(L(u_j))$. We denote by $Y$ the union of $X^{\sing}_{w_k}$ for
all $k$. Since the functions $u_j$ are bounded on $f^{-1}(W_2)$,
the argument of Proposition
\ref{P:ln} shows that
$$\il {Y}dd^cu_1\wedge\dots\wedge dd^cu_{n-1}\wedge dd^cv^*=0.$$
 Again
$\nu(Y)=0$.
\par Write $Z=X\cup Y$. It follows that
$$\nu(E)=-\il{E\sm Z}\phi\psi^* T\wedge dd^cv^*.$$
\par  We fix a point $w_0\in{\Bbb C}\sm W_1$ that we assume to be
equal to 0. Since $X_{0}\sm Z$ is a complex manifold, for every
its point there is a neighborhood and new coordinates
$\zeta=(\zeta_1,\dots,\zeta_n)$ on it such that
$f(\zeta)=\zeta_1^m$. A {\it preferred} neighborhood $U_{z}(\eps)$
of a point $z\in X_0\sm Z$ is an open set that in new
$\zeta$-coordinates has the form $\{|\zeta_1|<\eps^{1/m},
'\zeta=(\zeta_2,\dots,\zeta_n)\in {\Bbb D}^{n-1}(0,\eps)\}$.
\par Now, let us take an increasing sequence of nonnegative functions
$\phi_k\in\D^{0,0}(D\sm Z)$ converging to 1 on $D\sm Z$. There are
finitely many components $A_{kl}$, $1\le l\le l_k$, of $X_0$
intersecting the set $F_k=\supp\phi_k$.
\par Let us fix $k$ and choose a finite set of preferred  neighborhoods
$$U_{kli}=\{|\zeta_1|<\eps_{kli}^{1/m_{kl}},
'\zeta\in{\Bbb D}^{n-1}(0,\eps_{kli})\},$$ where $m_{kl}$ is the
multiplicity of $f$ on $A_{kl}$, covering each set $A_{kl}\cap
F_k$, and nonnegative functions $\psi_{kli}\in\D^{0,0}(U_{kli})$
such that $\sum_i\psi_{kli}\equiv1$ on a neighborhood $W_k$ of
$X_0\cap F_k$. Clearly there is an $\eps_0>0$ such that in the
preferred coordinates the set
$$U_{kli}(\eps)=\{\zeta\in U_{kli},|\zeta_1|<\eps^{1/m_{kl}},
'\zeta\in{\Bbb D}^{n-1}(0,\eps_{kli})\}\sbs W_k$$ when
$\eps\le\eps_0$. Note that $U_{kli}(\eps)=\{\zeta\in
U_{kli}:\,|f(\zeta)|<\eps\}$.
\par Let $V(\eps)=f^{-1}({\Bbb D}(0,\eps))$. There is
a positive $\eps_1<\eps_0$ such that
$$V(\eps)\cap F_k\sbs
U_k(\eps)=\bigcup\limits_{l,i}U_{kli}(\eps)$$ when $\eps<\eps_1$.
Hence, if the measure $\nu_k$ is defined by
$$\nu_k(E)=-\il{E}\phi\phi_k\psi^* T\wedge dd^cv^*,$$
then
\begin{equation}\label{e:sin}
\nu_k(V(\eps))=
-\sum_{l,i}\il{U_{kli}(\eps)}\phi\phi_k\psi_{kli}\psi^* T\wedge
dd^cv^*\end{equation} when $\eps<\eps_1$.
\par Since every preferred neighborhood $U_{kli}$ has
a structure of a direct product ${\Bbb
D}(0,\eps_{kli}^{1/m_{kl}})\times{\Bbb D}^{n-1}(0,\eps_{kli})$ and
$v^*(\zeta)=v(\zeta_1^{m_{kl}})=v_1(\zeta)$ on $U_{kli}$, by Lemma
\ref{L:dp}
$$\Psi_{kli}(\eps)=-\il{U_{kli}(\eps)}\phi\phi_k\psi_{kli}\psi^* T\wedge
dd^cv^*=-\il{{\Bbb D}(0,\eps^{1/m})}\psi
dd^cv_1\il{A_{\xi}}\phi\phi_k\psi_{kli}T,$$ where
$A_\xi=\{\zeta\in U_{kli}:\,\zeta_1=\xi\}$.
\par Since the functions $u_j$ are continuous and not equal to
$-\infty$ near $X_0$, the functions
$$\Phi_{kli}(\xi)=\il{A_\xi}\phi\phi_k\psi_{kli}T$$ are also
continuous.
\par If
$$\eta(E)=\il Edd^cv,$$ then
$$\il{{\Bbb D}(0,\eps^{1/m_{kl}})}dd^cv_1=
m_{kl}\eta({\Bbb D}(0,\eps)).$$ Hence,
$$\lim\limits_{\eps\to0}\frac{\Psi_{kli}(\eps)}{\eta({\Bbb
D}(0,\eps))}=m_{kl}\psi(0)\Phi_{kli}(0).$$ By (\ref{e:sin})
$$
\lim_{\eps\to0}\frac{\nu_k(V(\eps))}{\eta({\Bbb D}(0,\eps))}=
-\sum_{l,i}m_{kl}\psi(0)\Phi_{kli}(0).$$ But the last sum is equal
to
$$-\sum_Am_A\il{X_0}
\phi\phi_k\psi^*T=N_0(0,\phi\phi_k\psi^*)=N(0,\phi\phi_k\psi^*).$$
Thus, we see that the density of the measure
$\wtl\nu_k(E)=\nu_k(f^{-1}(E))$ with respect to $\eta$ is equal to
$N(0,\phi\phi_k\psi^*)$ at the points that do not belong to
$f(W_1)$. Since $v(w)$ is bounded near a point $w_0\in f(W_1)$,
its Laplacian has no mass at $w_0$. Hence,
$$\nu_k(D)=-\il{\Om}N(w,\phi\phi_k\psi^*)dd^cv.$$
\par The functions $\phi\phi_k\psi^*$ form an increasing sequence
and, therefore,
$$\lim_{k\to\infty}\nu_k(D)=
-\il{D}\phi\psi^* T\wedge dd^cv^*,$$ while
$$\lim_{k\to\infty}N_0(0,\phi\phi_k\psi^*)=N_0(0,\phi\psi^*)=
N(0,\phi\psi^*).$$ Hence
$$\il{D}\phi\psi^* T\wedge dd^cv^*=
\il{\Om}N(w,\phi\psi^*)dd^cv$$ and we proved the theorem in the
case of $\phi,\psi\ge0$ and $u_0\le0$ on $\supp\phi$.
\par For the proof of the general case we just note that the
function $u_0$ does not exceed some constant $a>0$ and, therefore,
is the difference of two plurisubharmonic functions negative on
$\supp\phi$. The functions $\phi$ and $\psi$ also can be written
as differences of nonnegative test functions. Hence, the general
case can be reduced to the case we have just considered.
\end{pf}

\section{Nevanlinna counting functions}\label{S:ncf}

\par Let $D\subset {\Bbb C}^n$ be a hyperconvex domain with a
continuous exhaustion function $u$ and let $f$ be a holomorphic
function on $D$. For $r<0$ we set $T_r=(u-r)(dd^cu)^{n-1}$. Let us
introduce the {\it Nevanlinna counting functions} as
$$N_{u,f}(w,r)=N(w,r)=-\il{B_u(r)}T_r\wedge dd^c\log|f-w|$$
and
$$N_0(w,r)=
-\sum_Am_{A,f}\il{A^{\reg}_{w,r}}T_r,$$ where the summation runs
over all irreducible components $A_w$ of $X_w$ and
$A^{\reg}_{w,r}=A^{\reg}_w\cap B_u(r)$.
\par Theorem \ref{T:cvi} leads to the following result.
\par
\bT\label{T:cvi1} If $w_0\not\in f(L(u))$ then
$N(w_0,r)=N_0(w_0,r)$ and the function $N(w,r)$ is lower
semicontinuous at $w_0$.
\par Let $\Om$ be a domain in ${\Bbb C}$ such that $f(D)\sbs\Om$
and let $v$ be a subharmonic function on $\Om$ and $v^*=v\circ f$.
If $L(u)$ is a finite set and $f(L(u))\cap L(v)=\emptyset$, then
$$\il {B_u(r)}(r-u)(dd^cu)^{n-1}\wedge dd^cv^*=
\il{\Om}N(w,r)dd^cv.$$
\eT
\begin{pf} In Theorem \ref{T:cvi} let $u_0=u-r$ and $u_1=\dots=u_{n-1}=u$.
Take increasing sequences of nonnegative functions
$\{\phi_k\}\sbs\D^{0,0}(D)$ with $\supp\phi_k\sbs B_u(r)$ and
$\{\psi_k\}\sbs\D^{0,0}({\Bbb D})$ converging to 1 on $B_u(r)$ and
$\Om$ respectively. Now Theorem \ref{T:cvi} and the monotone
convergence theorem yield the result.
\end{pf}
\par Let us explain why  we call these functions the Nevanlinna
counting functions. Our first observation is that by Fubini's
theorem
\begin{equation}\begin{align}\label{e:nct}
N(w,r)=&\il{B_u(r)}(r-u)(dd^cu)^{n-1}\wedge dd^c\log|f-w|
\notag\\&=\il{-\infty}^r\,dt\il{B_u(t)}(dd^cu)^{n-1}\wedge
dd^c\log|f-w|.\notag \end{align}\end{equation}
\par Therefore, if
$$n_{u,f}(w,r)=n(w,r)=\il {B_u(r)}(dd^cu)^{n-1}\wedge
dd^c\log|f-w|,$$ then
$$N(w,r)=\il{-\infty}^rn(w,t)\,dt.$$
\par  The function $n(w,t)$ has a nice geometric interpretation.
First of all, if $n=1$ then $n(w,r)$ is equal to the number of
points in $B_u(r)$, where $f(z)=w$, counted with their
multiplicities. Hence, it coincides with the function $n(w,r)$
from Nevanlinna's theory. If $D={\Bbb D}$ and $u(z)=\log|z|$, then
$$N(w,\log\rho)=\il{0}^\rho \frac{n(w,\log t)}t\,dt$$ and we see again
a classical formula.
\par If $n>1$ and $w\not\in f(L(u))$, then by Proposition
\ref{P:ln}
$$n(w,r)=\il{X_w\cap B_u(r)}(dd^cu)^{n-1}.$$
The form $(dd^cu)^{n-1}$ can be viewed as $2(n-1)$-dimensional
volume form on complex manifolds in the pseudo-metric $dd^cu$. For
example, if $D$ is the unit ball and $u(z)=|z|^2-1$, then
$(dd^cu)^{n-1}$ is exactly a scalar multiple of this volume form.
So in the multidimensional case, like in the classical parabolic
setting, we count not the number of preimages of $w$ in $B_u(r)$
but their pseudo-areas.
\par However, contrary to the classical parabolic case we
are more interested in the final form of the counting function,
i.e., the function
$$N_{u,f}(w)=N(w)=-\il Du(dd^cu)^{n-1}\wedge dd^c\log|f-w|,$$
which is the limit of $N(w,r)$ as $r\to0^-$.
\par For $\al>-1$ and $u<0$ we introduce the following auxiliary functions
$$\si_\al(u)=\il{u}^0|r|^\al e^r\,dr=
\frac1{\al+1}|u|^{\al+1}+o(|u|^{\al+1})$$ and
$$\gm_\al(u)=\il{u}^0|r|^\al e^r(r-u)\,dr=
-\si_{\al+1}(u)-u\si_\al(u).$$ Note that both functions are
positive,  $\si_\al(u)\le\Gm(\al+1)$ and
$$\gm_\al(u)=\frac1{(\al+1)(\al+2)}|u|^{\al+2}+o(|u|^{\al+2}).$$
We also set $\si_{-1}(u)=1$ and $\gm_{-1}(u)=-u$.
\par We define  the Nevanlinna counting function of order $\al$ as
$$N_{u,f,\al}(w)=N_\al(w)=
\il {D}\gm_\al(u)(dd^cu)^{n-1}\wedge dd^c\log|f-w|.$$ By Fubini's
theorem
\begin{equation}\label{e:dna}
N_\al(w)=\il{-\infty}^0|r|^\al e^rN(w,r)\,dr.
\end{equation}
\par If $n=1$ then $N_\al(w)=\sum_{f(z)=w}\gm_\al(u(z))$.
\par The importance of the Nevanlinna counting function is
explained by the fact that the Hardy and weighted Bergman norms
of functions can be expressed in terms of their Nevanlinna
counting functions. The following result may be viewed as
a generalization of the classical Littlewood-Paley formula.
\bT\label{T:cvn} Let $D$ be a hyperconvex domain in ${\Bbb C}^n$
with an exhaustion function $u$ such that the set $L(u)$ is
finite. If $f$ is a holomorphic function on $D$, then
\begin{equation}\begin{align}
\|f\|^p_{A^p_{u,\al}}-\il D\si_\al(u)|f|^p(dd^cu)^n&=
\il{D}\gm_\al(u)dd^c|f|^p\wedge(dd^cu)^{n-1}\notag\\&=\il{\Bbb
C}N_\al(w)dd^c|w|^p. \notag
\end{align}\end{equation}
\eT
\begin{pf} The case $\al=-1$ follows immediately from the
definitions of the Hardy norm and the function $N(w)$ and from
Theorem \ref{T:cvi} where $v(w)=|w|^p$.
\par When $\al>-1$ we recall that
\begin{equation}\begin{align}
&\|f\|^p_{A^p_{u,\al}}=\il{-\infty}^0|r|^\al
e^r\mu_{u,r}(|f|^p)\,dr= \il{-\infty}^0|r|^\al
e^r\left(\il{B_u(r)} |f|^p(dd^cu)^n\right)\,dr\notag\\&+
\il{-\infty}^0|r|^\al
e^r\left(\il{B_u(r)}(r-u)dd^c|f|^p\wedge(dd^cu)^{n-1}\right)\,dr.
\notag
\end{align}\end{equation}
By Fubini's theorem the first double integral is equal to
$$\il D\si_\al(u)|f|^p(dd^cu)^n$$ and the second integral
is equal to
$$\il{D}\gm_\al(u)dd^c|f|^p\wedge(dd^cu)^{n-1}.$$
By Theorem \ref{T:cvi} the second integral is equal to
$$ \begin{array}{c} \il{-\infty}^0|r|^\al
e^r\left(\il{\Bbb C}N(w,r)dd^c|w|^p\right)\,dr=\il{\Bbb
C}\left(\il{-\infty}^0|r|^\al e^rN(w,r)\,dr\right)\,dd^c|w|^p \\
   \\   =\il{\Bbb C}
N_\alpha (w)dd^c|w|^p. \end{array} $$

\end{pf}
\par In the case when $u=g_D(z,z_0)$ is a pluricomplex Green
function the theorem above takes exactly the form of the
Littlewood--Paley identity.
\bC\label{C:lpt} If in the assumptions of Theorem \ref{T:cvn} the
function $u=g_D(z,z_0)$, then
$$\|f\|^p_{A^p_{u,\al}}=\Gm(\al+1)|f(z_0)|^p+
\il{\Bbb C}N_\al(w)dd^c|w|^p$$
\eC

\section{Properties of the Nevanlinna counting functions}\label{S:pncf}

\par By their definition the Nevanlinna counting functions depend
on the choice of the exhaustion. The next result shows that this
dependence can be estimated.
\bT\label{T:enct}Let $u$ and $v$ be two continuous exhausting
functions on $D$ such that $cv\le u\le c^{-1}v$ for all $z\in D$.
If $\phi$ is a plurisubharmonic function on $D$ and
$$\il{B_u(r)}(r-u)dd^c\phi\wedge(dd^c(u+v))^{n-1}<\infty$$ for
all $r<0$, then for all $r<0$
$$\il{B_u(r)}(r-u) dd^c\phi\wedge(dd^cu)^{n-1}\le
c^n\il{B_v(r_n)}(r_n-v) dd^c\phi\wedge(dd^cv)^{n-1},$$ where
$r_n=c^{-2n+1}r$.
\eT
\begin{pf} The inequality $cv\le u\le c^{-1}v$ on $D$ implies two
facts: first, $B_v(cr)\sbs B_u(r)\sbs B_v(r/c)$ and, secondly,
$c(c^{-1}r-v)\ge r-u$. Thus
\begin{equation}\label{e:s0}
\il{B_u(r)}(r-u)dd^c\phi\wedge(dd^cu)^{n-1}\le
c\il{B_v(c^{-1}r)}(c^{-1}r-v)dd^c\phi \wedge(dd^cu)^{n-1}.
\end{equation}

\par Since
$(dd^c(u+v))^{n-1}\ge (dd^cv)^k\wedge (dd^cu)^{n-k-1}$, $0\le k\le
n-1$, we conclude from the assumptions of the theorem that
\begin{equation}\label{e:fs}
\il {B_u(r)}dd^c\phi\wedge (dd^cv)^k\wedge (dd^cu)^{n-k-1}<\infty
\end{equation}
for all $r<0$.
\par If $T= dd^c\phi\wedge(dd^cv)^k\wedge(dd^cu)^{n-k-2}$, then
$T$ is a closed positive current of bidimension $(1,1)$ and by
(\ref{e:fs})
$$\il {B_v(r)}dd^cu\wedge T<\infty.$$
\par Since the function $c^{-1}r-u\ge0$ on $B_v(r)$, we have by
Theorem \ref{T:C}
$$\il {B_v(r)}(r-v)(dd^cu)\wedge T\le
\il {B_v(r)}(c^{-1}r-u)(dd^cv)\wedge T.$$ But $c^{-1}r-u\le
c(c^{-2}r-v)$ and, therefore, the latter integral does not exceed
$$c\il{B_v(c^{-2}r)}(c^{-2}r-v)dd^c\phi\wedge(dd^cv)^{k+1}
\wedge(dd^cu)^{n-k-2}.$$ Thus,
\begin{equation}\label{e:ie}\begin{align}
&\il {B_v(r)}(r-v)\wedge dd^c\phi\wedge
(dd^cv)^k\wedge(dd^cu)^{n-k-1}\notag\\&\le
c\il{B_v(r_1)}(r_1-v)dd^c\phi\wedge(dd^cv)^{k+1}
\wedge(dd^cu)^{n-k-2},
\end{align}\end{equation} where
$r_1=c^{-2}r$.
\par Starting with (\ref{e:s0}) and applying (\ref{e:ie}) $n-1$
times to the right-hand side of (\ref{e:s0}) to eliminate $dd^cu$
we get the theorem.
\end{pf}

\par Since by  Theorem \ref{P:ln} $N_{u+v,f}(w,r)<\infty$ when
$L(u)\cap X_w=\emptyset$ and $L(v)\cap X_w=\emptyset$, in the case
when $\phi(z)=\log|f(z)-w|$, the previous result   and
(\ref{e:dna})
 imply the following corollary.

\bC\label{C:inct} When
$L(u)\cap X_w=\emptyset$ and $L(v)\cap X_w=\emptyset$, in the
assumptions of Theorem \ref{T:enct} $c^{-n}N_v(w,r_n)\le N_u(w,r)\le c^nN_v(w,r_n)$,
where $r_n=c^{-2n+1}r$. Consequently, $c^{-n}N_{v,\al}(w)\le N_{u,\al}(w)\le
c^nN_{v,\al}(w)$
\eC
\par It is well-known that for a self-mapping of the unit disk $f$ which
fixes the origin, the classical Nevanlinna counting function satisfies
the estimate $N(w)\leq -\log|w|$. Our next result generalizes this
statement to the multivariable case.

\bT\label{T:mnct} Let $u\in\E(D)$ and let $f:D\to{\Bbb D}$ be a
holomorphic function. Then there is a number $r<1$ and a constant
$c$ depending only on $u$ and $f$ such that for every $w\in{\Bbb
D}$, $|w|>r$, $N_u(w)\le -c\log|w|$. \par Moreover if
$u(z)=g_D(z,z_0)$ and $f(z_0)=w_0$, then
$$N_u(w)\le-(2\pi)^n\log\left|\frac{w_0-w}{1-\ovr
w_0w}\right|.$$
\eT
\begin{pf} We fix a point $z_0$ in $D$ and let $v(z)=g_D(z,z_0)$.
If $w_0=f(z_0)$ we set $v_1(\zeta)=\log|\zeta-w_0|/|1-\ovr
w_0\zeta|$ and $\phi(\zeta)=\log|\zeta-w|/|1-\ovr w\zeta|, \
\zeta\in {\Bbb D}$. By Theorem \ref{T:lmp} $\mu_{v,r}(\phi\circ
f)\le(2\pi)^{n-1}\mu_{v_1,r}(\phi)$. If $w\ne w_0$, then
$$\mu_{v,r}(\phi\circ f)=(2\pi)^n\log\left|\frac{w_0-w}{1-\ovr
w_0w}\right|+N_v(w,r)$$ and
$$\mu_{v_1,r}(\phi)=2\pi\log\left|\frac{w_0-w}{1-\ovr
w_0w}\right|+\il{B_{v_1}(r)}(r-v_1)dd^c\phi=2\pi r.$$
Hence,
$$N_v(w,r)\le(2\pi)^n r-(2\pi)^n\log\left|\frac{w_0-w}{1-\ovr
w_0w}\right|$$ and
$$N_v(w)\le-(2\pi)^n\log\left|\frac{w_0-w}{1-\ovr
w_0w}\right|.$$
\par For an arbitrary function $u\in\E(D)$ we take the compact
set $K=\ovr B_v(-1)\cup\ovr{\{u<-1\}}$. Let $u_1=\max\{u,-1\}$ and
$v_1=\max\{v,-1\}$. If $f(K)\sbs{\Bbb D}(0,r)$ and $|w|>r$, then
$N_u(w,r)=N_{u_1}(w,r)$ and $N_v(w,r)=N_{v_1}(w,r)$. Since there
is a constant $c>1$ such that $cv_1\le u_1\le c^{-1}v_1$ for all
$z\in D$, by Corollary \ref{C:inct}
$$N_u(w,r)\le c^nN_v(w,r_n)\le(2\pi c)^nr_n-(2\pi c)^n\log\left|\frac{w_0-w}{1-\ovr
w_0w}\right|,$$ where $r_n=c^{-2n+1}r$. Hence,
$$N_u(w)\le-(2\pi c)^n\log\left|\frac{w_0-w}{1-\ovr
w_0w}\right|\le-(2\pi c_1)^n\log|w|$$ when $|w|>r$.
\end{pf}
\par Our next result establishes a multidimensional analog of Shapiro's mean value inequality
for counting functions.
\bT\label{T:smi} Let $f$ be a analytic function on a hyperconvex
domain $D\sbs{\Bbb C}^n$ with an exhausting function $u$, and let
$K=\supp(dd^cu)^n$ be a compact set. If $w_0\in{\Bbb C}$ and
$0<\rho<r_0=\dist(w_0,f(K))$, then
$$N_{u,\al}(w_0)\le\frac1{2\pi i\rho^2}\il{{\Bbb D}(w_0,\rho)}
N_{u,\al}(w)\,dw\wedge d\ovr w.$$
\eT
\begin{pf} Note that the function
$\Phi(w,r)=\mu_{u,r}(\log|f(z)-w|)$ is subharmonic and
$$\Phi(w_0,r)\le\frac1{2\pi it^2}
\il{{\Bbb D}(w_0,t)}\Phi(w,r)\,dw\wedge d\ovr w.$$ The function
$$\Psi(w)=\il D\log|f(z)-w|(dd^cu)^n$$ is harmonic outside of $f(K)$.
By (\ref{e:ljf}) $N_{u}(w,r)=\Phi(w,r)-\Psi(w)$ when $K\sbs
B_u(r)$. Hence, the function $N_{u}(w,r)$ is subharmonic outside of
$f(K)$. Since this family of functions is increasing in $r$ and
$N_{u}(w)=\lim_{r\to0^-}N_{u}(w,r)$, the theorem follows from
(\ref{e:dna}).
\end{pf}
\par This theorem has an important corollary.
\bC\label{C:smi} Let $f$ be an analytic function on a hyperconvex
domain $D\sbs{\Bbb C}^n$ with an exhausting function $u\in\E(D)$,
and let $K\sbs D$ be a compact set such that $L(u)$ lies in the
interior of $K$. There is a constant $c>0$ depending only on $u$
and $F$ such that if $w_0\in{\Bbb C}$ and
$0<\rho<r_0=\dist(w_0,f(K))$, then
$$N_{u,\al}(w_0)\le\frac{c}{2\pi i\rho^2}\il{{\Bbb D}(w_0,\rho)}
N_{u,\al}(w)\,dw\wedge d\ovr w.$$
\eC
\begin{pf} Let $z_0$ be an interior point of $K$. We take the
function $g_D(z,z_0)$, let $a$ be the maximum of $g_D(z,z_0)$ on a
closed ball $B$ lying in the interior of $F$ and set
$v(z)=\max\{g_D(z,z_0),a\}$. Then $\supp(dd^cv)^n$ lies in the
interior of $K$.
\par If $b$ is the maximum of $g_D(z,z_0)$ on $B$ and
$u_1(z)=\max\{u(z),b\}$, then there is a constant $k\ge1$ such
that $k^{-1}v\le u_1\le kv$ on $D$. By Corollary \ref{C:inct}
$k^{-n}N_{u_1,\al}\le N_{v,\al}\le k^nN_{u_1,\al}$. Since
$N_{u_1,\al}(w)=N_{u,\al}(w)$ when $w\not\in f(K)$ and
$$N_{v,\al}(w_0)\le\frac{c}{2\pi i\rho^2}\il{{\Bbb D}(w_0,\rho)}
N_{v,\al}(w)\,dw\wedge d\ovr w,$$ the corollary follows.
\end{pf}

\section{Application to composition operators: boundedness and compactness}\label{S:cobc}

\par Now, let $F:D_1\to D_2$ be a holomorphic mapping between
hyperconvex domains $D_1\subset{\Bbb C}^n$ and $D_2\subset{\Bbb
C}^m$ with exhausting functions $u_1\in {\E}(D_1)$ and $u_2\in
{\E}(D_2)$ respectively. If $f$ is a holomorphic function on $D_2$
then we denote by $f^*$ the function $f\circ F$.
\par Let us introduce the ``tail'' part of the Nevanlinna function
$$N^*_{u_1,F,f,\beta}(w,r)=
\il {T(r)}\gm_\beta(u_1)(dd^cu_1)^{n-1}\wedge dd^c\log|f^*-w|,$$
where $T(r)=D_1\sm \ovr B_{u^*_2}(r)=\{z\in D_1:\,u_2(F(z))>r\}$. We
define the $(\beta, \alpha)$-{\em deficiency} of $F$ as
$$\dl_{u_1,u_2,F,\beta,\al}(r)=\dl_{F,\beta,\al}(r)=
\sup\frac{N^*_{u_1,F,f,\beta}(w,r)}{N_{u_2,f,\al}(w)},$$ where the
supremum is taken over all $f\in A^p_{u_2,\al}(D_2)$ and all
$w\in{\Bbb C}$. We assume that the ratio
$N^*_{u_1,F,f,\beta}(w,r)/N_{u_2,f,\al}(w)=0$ when
$N_{u_2,f,\al}(w)=0$. Clearly, the function $\dl_{F,\al,\beta}(r)$
is decreasing in $r$.
\par In what follows, given a holomorphic function $h$ on a domain $D$
and coordinates $\zeta_1,...\zeta_n$ in $\aC ^n$ we denote by
$\nabla h$ the complex gradient of $h$,
$$\nabla h=\left(\frac{\bd h}{\bd\zeta_1},\dots,
\frac{\bd h}{\bd\zeta_n}\right),$$ and for vectors
$a=(a_1,\dots,a_n)$ and $b=(b_1,\dots, b_n)$ in ${\Bbb C}^n$ we
let $\langle a,b\rangle=\sum a_ib_i$.

\bL\label{L:dfe} Let $F$ be a holomorphic mapping from a domain
$D_1\sbs{\Bbb C}^n$ into  domain $D_2\sbs{\Bbb C}^m$, and let $h$
be a holomorphic function on $D_2$. If $h^*=h\circ F$ and
$F=(F_1,\dots,F_m)$, then at every point $\zeta\in D_1$ the
following estimate holds
$$dd^c|h^*|^p\le\frac{p^2}4|h^*|^{p-2}(\zeta)|\nabla
h(F(\zeta))|^2\sum_{i=1}^mdd^c|F_i(\zeta)|^2.$$
\eL
\begin{pf} A direct calculation for a holomorphic function $h$
shows that
$$dd^c|h|^p=\frac{p^2}4|h|^{p-2}dd^c|h|^2$$
and  $dd^c|h^*|^2=2i\bd h^*\wedge\ovr\bd \ovr h^*$. The quadratic form
corresponding to $dd^c|h^*|^2$ is
$$\sum_{i,j=1}^n\frac{\bd^2|h^*|^2}{\bd\zeta_i\bd\ovr\zeta_j}(\zeta)\xi_i\ovr\xi_j=\sum_{i,j=1}^n\frac{\bd
h^*}{\bd\zeta_i}(\zeta) \frac{\bd\ovr
h^*}{\bd\ovr\zeta_j}(\zeta)\xi_i\ovr\xi_j=|\langle\nabla
h^*(\zeta),\xi\rangle|^2,$$where $\xi=(\xi_1,\dots,\xi_n)$. Now
$$|\langle\nabla h^*(\zeta),\xi\rangle|=|\langle\nabla
h,F'(\zeta)\xi\rangle|\le|\nabla h(F(\zeta))||F'(\zeta)\xi|.$$ But
$$|F'(\zeta)\xi|^2=\sum_{i=1}^m|\langle\nabla
F_i(\zeta),\xi\rangle|^2.$$ Thus,
$$|\langle\nabla h^*(\zeta),\xi\rangle|^2\le|\nabla
h(F(\zeta))|^2\sum_{i=1}^m|\langle\nabla F_i(\zeta),\xi\rangle|^2.$$
\par Since
$$dd^ch^*=2i\sum_{i,j=1}^n\frac{\bd h^*}{\bd\zeta_i}
\frac{\bd\ovr h^*}{\bd\ovr\zeta_j}d\zeta_i\wedge d\ovr\zeta_j$$
and
$$dd^c|F_i|^2=2i\sum_{i,j=1}^n\frac{\bd F_i}{\bd\zeta_i}
\frac{\bd\ovr F_i}{\bd\ovr\zeta_j}d\zeta_i\wedge d\ovr\zeta_j,$$
by Corollary 3.2.5 from (\cite{Kl}, p.102) we conclude that the
differential form
$$|\nabla
h(F(\zeta))|^2\sum_{i=1}^mdd^c|F_i(\zeta)|^2-dd^c|h^*(\zeta)|^2$$
is positive.
\end{pf}
\bL\label{L:ie10} Let $F:D_1\to D_2$ be a holomorphic mapping from
a hyperconvex domain $D_1\subset{\Bbb C}^n$ with an exhausting
functions $u$ such that the set $L(u)$ is finite, into a domain
$D_2\subset{\Bbb C}^m$. If $V$ is an open set in $D_1$, $\phi$ is
a nonnegative continuous function on $D_1$, $W$ is an open set in
$D_2$, $F(V)\sbs W$ and $f$ is a holomorphic function on $W$ such
that $0<c_1<|f(z)|<c_2$ and $|\nabla f|<c_3$ on $W$, then
$$\il{V}\gm_\beta(u)dd^c|f^*|^p\wedge(dd^cu)^{n-1}\le
c^{p-2}c_3^2\Gm(\beta+1)\sum_{i=1}^m\|F_i\|^2_{H^2_{u}(D_1)},$$
where the constant $c$ is equal to $c_2$ when $p-2\ge0$ and $c_1$
when $1\le p\le 2$.
\begin{pf}\par The integral in question is equal to
$$\frac{p^2}4\il{V}\gm_\beta(u)
|f^*|^{p-2}dd^c|f^*|^2\wedge(dd^cu)^{n-1}.$$ Since
$0<c_1<|f(z)|<c_2$ and $|\nabla f|<c_3$ on $V$, by Lemma
\ref{L:dfe} this integral does not exceed
$$c^{p-2}c_3^2\sum_{i=1}^m\il{D_1}
\gm_\beta(u)dd^c|F_i|^2\wedge(dd^cu)^{n-1},$$
\par By Theorem \ref{T:cvn}
$$\il{D_1}\gm_\beta(u)dd^c|F_i|^2\wedge(dd^cu)^{n-1}=
\il{-\infty}^0|r|^\beta
e^r\il{B_{u}(r)}(r-u)dd^c|F_i|^2\wedge(dd^cu)^{n-1}.$$ Since
$$\il{B_{u}(r)}(r-u)dd^c|F_i|^2\wedge(dd^cu)^{n-1}\le
\|F_i\|^2_{H^2_{u}(D_1)},$$ we see that
$$\il{V}\gm_\beta(u)dd^c|f^*|^p\wedge(dd^cu)^{n-1}\le
c^{p-2}c_3^2\Gm(\beta+1)\sum_{i=1}^m\|F_i\|^2_{H^2_{u}(D_1)}.$$
\end{pf}
\eL
\par The following result gives sufficient conditions of boundedness
and compactness of the composition operator $C_F$ as an operator
acting from $A^p_\alpha (D_2)$ into $A^p_\beta (D_1)$. It could be
viewed as an extension of the appropriate results of Shapiro
\cite{Sha2} and Smith \cite{Sm} to the multi-variable case.

\bT\label{T:bct} Let $F:D_1\to D_2$ be a holomorphic mapping between
hyperconvex domains $D_1\subset{\Bbb C}^n$ and $D_2\subset{\Bbb
C}^m$ with exhausting functions $u_1$ and $u_2$ respectively such
that the sets $L(u_1)$ and $L(u_2)$ are finite. \be \item If there
exists $r_0<0$ such that $ \dl_{F,\al,\beta}(r_0)<\infty $, then
$C_F$ is a bounded operator from $A^p_{u_2,\alpha}(D_2)$ into
$A^p_{u_1,\beta}(D_1)$. \item  If the function
$\dl_{F,\al,\beta}(r)$ converges to $0$ as $r\to0$, then $C_F$ is
a compact operator from $A^p_{u_2,\alpha}(D_2)$ into
$A^p_{u_1,\beta}(D_1)$. \ee\eT
\begin{pf} To prove the first part of the theorem we take a
function $f\in A^p_\al(D_2)$ with $\|f\|_{A^p_\al(D_2)}=1$. If $K$
is a compact set in $D_1$ such that $(dd^cu_1)^n\equiv0$ outside
of $K$ and $K'=F(K)$, then by Theorem \ref{T:ci} there is a
constant $C_1$ not depending on $f$ such that $|f(z)|\le C_1$ on
$K'$. Consequently,
\begin{equation}\label{e:bct1}
\il {D_1}\si_\beta(u_1)|f^*|^p(dd^cu_1)^n\le C_1^p\Gm(\beta+1)\il
{D_1}(dd^cu_1)^n.
\end{equation}
\par By Theorem \ref{T:ci} there are constants $C_2, C_3>0$ not
depending on $f$ such that $|f(z)|\le C_2$ and $|\nabla f(z)|^2\le
C_3$ when $z\in B_{u_2}(r_0/2)$. Let $g(z)=f(z)+2C_2$. Then
$C_2<|g(z)|<3C_2$ on $B_{u_2}(r_0/2)$ and
$$\|g\|_{A^p_\al(D_2)}\le
1+(2C_2)\left ( \il {D_2}\si_\al(u_2)(dd^cu_2)^n\right )
^{1/p}=C_4.$$
\par By Lemma \ref{L:ie10} we have
$$\il{B_{u^*_2}(r_0/2)}
\gm_\beta(u_1)dd^c|g^*|^p\wedge(dd^cu_1)^{n-1}\le
c^{p-2}C_3^2\Gm(\beta+1)\sum_{i=1}^m\|F_i\|^2_{H^2_{u_1}(D_1)}.$$
But the functions $F_i$ are bounded and, therefore,
$\|F_i\|^2_{H^2_{u_1}(D_1)}<\infty$. Hence
\begin{equation}\label{e:bct3}
\il{B_{u^*_2}(r_0/2)}
\gm_\beta(u_1)dd^c|g^*|^p\wedge(dd^cu_1)^{n-1}<C_5,
\end{equation}
where the constant $C_5$ depends only on $c, C_3,p$ and $F$.
\par By Theorem \ref{T:cvi}
$$\il{T(r_0)}
\gm_\beta(u_1)dd^c|g^*|^p\wedge(dd^cu_1)^{n-1}=\il{\Bbb
C}N^*_{u_1,F,g,\beta}(w,r_0)dd^c|w|^p.$$ Since
$N^*_{u_1,F,g,\beta}(w,r_0)\le
\delta_{u_1,u_2,F,\beta,\alpha}(r_0)N_{u_2,g,\al}(w)$,
\begin{eqnarray}
\il{\Bbb C}N^*_{u_1,F,f,\beta}(w,r_0)dd^c|w|^p\le
\delta_{u_1,u_2,F,\beta,\alpha}(r_0) \il{\Bbb
C}N_{u_2,g,\al}(w)dd^c|w|^p \nonumber \\
\le
\delta_{u_1,u_2,F,\beta,\alpha}(r_0)\|g\|^p_{A^p_{u_2,\al}(D_2)}.
\label{e:bct6}
\end{eqnarray}
Combining together (\ref{e:bct1}) and (\ref{e:bct3}) we see that
$$\|g^*\|^p_{A^p_{u_1,\beta}(D_1)}\le
\delta_{u_1,u_2,F,\beta,\alpha}(r_0)\|g\|^p_{A^p_{u_2,\al}(D_2)}.$$ Since
$$\|f^*\|^p_{A^p_{u_1,\beta}(D_1)}\le
\|g^*\|^p_{A^p_{u_1,\beta}(D_1)}+2C_2\|1\|^p_{A^p_{u_1,\beta}(D_1)},$$ this
implies that the composition operator is bounded.
\par  Suppose that the second condition is satisfied. To show the
compactness of $C_f$ we need to show that if holomorphic functions
$g_k$ converge to 0 uniformly on compacta in $D_2$ and
$\|g_k\|_{A^p_\al(D_2)}\le1$, then
$\|g^*_k\|_{A^p_\beta(D_1)}\to0$ as $k\to\infty$.
\par For $\eps>0$ we take $r_0<0$ such that
$\dl_{F,\al,\beta}(r_0)<\eps$ and let $r=r_0/2$. We also take
$r_1>r$ such that the set $B_{u_1}(r)$ compactly belongs to
$B_{u_1}(r_1)$.
\par Let $\eps_k(r_1)$ be the supremum of $|g_k|$ on
$B_{u_2}(r_1)$. Then by Cauchy inequalities there is a constant
$C_1$ such that $|\nabla g_k|<C_1\eps_k(r_1)$ on $B_{u_2}(r)$. If
$h_k=2\eps_k(r_1)+g_k$ then $0<\eps_k(r_1)<|h_k|<3\eps_k(r_1)$ and
$|\nabla h_k|<C_1\eps_k(r_1)$ on $B_{u_2}(r)$. By Lemma
\ref{L:ie10}
$$\il{B_{u^*_2}(r)}
\gm_\beta(u_1)dd^c|h_k^*|^p\wedge(dd^cu_1)^{n-1}\le
Cc^{p-2}C_1^2\eps^2_k(r_1),$$ where
$C=\Gm(\beta+1)\sum_{i=1}^m\|F_i\|^2_{H^2_{u_1}(D_1)}$ and the
constant $c$ is equal to $3\eps_k(r_1)$ when $p-2\ge0$ and
$\eps_k(r_1)$ when $1\le p\le 2$. In any case the integral does
not exceed $9CC_1^2\eps^p_k(r_1)$. Thus
\begin{equation}\label{e:bct5}
\lim_{k\to\infty}\il{B_{u^*_2}(r)}
\gm_\beta(u_1)dd^c|h_k^*|^p\wedge(dd^cu_1)^{n-1}=0.
\end{equation}
By (\ref{e:bct6}) $$\il{T(r_0)}
\gm_\beta(u_1)dd^c|h_k^*|^p\wedge(dd^cu_1)^{n-1}\le
\dl_{F,\al,\beta}(r_0)\|h_k\|^p_{A^p_{u_2,\al}(D_2)}<a_k\eps,
$$
where $a_k=\|\eps_k(r_1)\|^p_{A^p_{u_2,\al}(D_2)}$. Combining
(\ref{e:bct5}) with the latter estimate we get that
$$\limsup_{k\to\infty}\il{D_1}
\gm_\beta(u_1)dd^c|h_k^*|^p\wedge(dd^cu_1)^{n-1}\le\eps.$$
\par Evidently,
$$\lim_{k\to\infty}\il {D_1}\si_\beta(u_1)|g^*_k|^p(dd^cu_1)^n=0,$$
and, therefore, the functions $h_k$ and, consequently, by Theorem
\ref{T:cvn} $g^\ast_k$, converge to 0 in $A^p_{u_1,\beta}(D_1)$
\end{pf}
\par To provide necessary conditions we fix a compact set
$K\sbs D_1$ whose interior contains $L(u_1)$ and for a holomorphic
function $f\in A^p_{u_2,\al}(D_2)$ introduce the function
$$\nu_{F,\al,\beta}(w,f)=\nu(w,f)=\frac{|w|^pN_{u_1,f^*,\beta}(w)}
{\|f\|^p_{A^p_{u_2,\al}}}$$ (here and below we use the same
notation: for a function $h$ on $D_2$ we denote by $h^\ast $ the
composition $h$ and $F$, \ $h^\ast=h\circ F$). For $a>1$ we set
$$\rho_{u_1,u_2,F,\al,\beta}(a)=\rho(a)=\sup\nu_{F,\al,\beta}(w,f),$$
where the supremum is taken over all $f\in A^p_{u_2,\al}(D_2)$ and
all $w\in{\Bbb C}$, $|w|>a\max_{\zeta\in K}|f^*(\zeta)|$. Note
that $\nu_{\al,\beta}(cw,cf)=\nu_{\al,\beta}(w,f)$. Thus, we may
assume that $\|f\|_{A^p_{u_2,\al}(D_2)}=1$ in the definition of
$\rho_{F,\al,\beta}$.
\par We also will need another characteristic of the mapping $F$.
For its definition we fix an open ball $B\sbs\sbs D_2$ and for $t>0$ and
$a>1$ we set
$$\wtl\rho_{u_1,u_2,F,\al,\beta}(t,a)=\wtl\rho_{F,\al,\beta}(t,a)=
\sup\nu_{\al,\beta}(w,f),$$ where the supremum is taken over all
$w\in{\Bbb C}$, $|w|>a\max_{\zeta\in K}|f^*(\zeta)|$, and all
$f\in A^p_{u_2,\al}(D_2)$ such that $\|f\|_{A^p_{u_2,\al}(D_2)}=1$
and $|f|<t$ on $B$.

\bT\label{T:bctn} Let $F:D_1\to D_2$ be a holomorphic mapping between
hyperconvex domains $D_1\subset{\Bbb C}^n$ and $D_2\subset{\Bbb
C}^m$ with exhausting functions $u_1(\zeta)\in\E(D_1)$ and $u_2\in
{\E}(D_2)$ respectively and such that the sets $L(u_1)$ and
$L(u_2)$ are finite. \be \item If $C_F$ is a bounded operator from
$A^p_{u_2,\alpha}(D_2)$ into $A^p_{u_1,\beta}(D_1)$, then
$\rho_{F,\al,\beta}(a)<\infty$ for all $a>1$. \item  If $C_F$ is a
compact operator from $A^p_{u_2,\alpha}(D_2)$ into
$A^p_{u_1,\beta}(D_1)$, then the function
$\wtl\rho_{F,\al,\beta}(t)$ converges to $0$ as $t\to0^+$ for all
$a>1$. \ee\eT
\begin{pf} If the first part of the theorem does not hold, then
there are functions $f_j\in A^p_{u_2,\al}(D_2)$,
$\|f_j\|_{A^p_{u_2,\al}}=1$, and $w_j\in{\Bbb C}$ such that
$|w_j|>a\max_{\zeta\in K}|f_j^*(\zeta)|$ and
$$|w_j|^pN_{u_1,f_j^*,\beta}(w_j)\ge
j.$$
\par Let $r_j=C_0|w_j|$, where $C_0=(1-a^{-1})/2$ and
$V_j={\Bbb D}(w_j,r_j)\subset \aC$.  Since $|w|\ge C_1|w_j|$ on
$V_j$, where $C_1=(a^{-1}+1)/2$, we get
\begin{equation}\begin{align} &\|f^*_j\|^p_{A^p_{u_1,\beta}}\ge
\il{V_j}N_{u_1,f^*_j,\beta}(w)\,dd^c|w|^p=\frac{p^2}4
\il{V_j}|w|^{p-2}N_{u_1,f^*_j,\beta}(w)\,dd^c|w|^2\notag\\&\ge
\frac{p^2}4C_1^p|w_j|^p
\il{V_j}|w|^{-2}N_{u_1,f^*_j,\beta}(w)\,dd^c|w|^2.\notag
\end{align}\end{equation}
But $|w|\le C_2|w_j|$ on $V_j$, where $C_2=(3+a^{-1})/2$. Hence,
$$\|f^*_j\|^p_{A^p_{u_1,\beta}}\ge \frac{p^2\pi
C_1^pC_0^2|w_j|^p}{2C_2^2}\frac1{2\pi r_j^2}
\il{V_j}N_{u_1,f^*_j,\beta}(w)\,dd^c|w|^2.$$ But the disk $V_j$
lies outside of $f_j^*(K)$. Thus, by Corollary \ref{C:smi}
$$\|f^*_j\|^p_{A^p_{u_1,\beta}}\ge\frac{\pi p^2C_1^pC_0^p}{2C_2^2}|w_j|^pN_{u_1,f^\ast_j,\beta}(w_j)
 \ge Cj,$$ where the constant $C=\frac{\pi p^2C_1^pC_0^p}{2C_2^2}$
depends only on $a$, $K$, $\lambda $ and $p$. Therefore, the norms
$\|f^*_j\|_{A^p_{u_1,\beta}}\to\infty$ as $j\to\infty$ and we get
a contradiction.
\par For the proof of the second part we assume that the statement
does not hold and take sequences of complex numbers $\{w_j\}$ and
functions $f_j\in A^p_{u_2,\al}$, $\|f\|^p_{A^p_{u_2,\al}}=1$,
such that $|w_j|>a\max_{\zeta\in K}|f_j^*(\zeta)|$, $|f_j|<1/j$ on
$B$ and $\nu_{\al,\beta}(w_j,f_j)\ge c>0$. Clearly, the sequence
$f_j$ converges to 0 uniformly  on compacta in $D_2$. But the same
estimates as in the first part show that
$\|f^*_j\|_{A^p_{u_1,\beta}}\ge C>0$ and this contradicts the
compactness of $C_F$.
\end{pf}

\section{Composition operators induced by  mappings into the unit
disk}\label{S:coimud}

\par It turns out that if $D_2$ is the unit disk, the necessary and
sufficient conditions, given by Theorems \ref{T:bctn} and
\ref{T:bct} respectively, agree. Of course, this is caused by the
simple structure of divisors in $\aC$. In this case the results of
the previous section allow us to state and prove a theorem which
gives necessary and sufficient conditions for boundedness and
compactness of a composition operator in the form of results of
\cite{Sha2} and \cite{Sm}.

\bT\label{T:bctd} Let $F:D\to {\Bbb D}$ be a holomorphic mapping
from a hyperconvex domains $D\subset{\Bbb C}^n$ into the unit disk
${\Bbb D}$ with exhausting functions $u\in {\E}(D)$ such that the
set $L(u)$ is finite and $v(z)=\log|z|$ respectively. Then the
condition
$$N_{F,u,\beta}(z)=O(\gm_\al(\log|z|))\text{\qquad as\qquad}
|z|\to1$$ is necessary and sufficient for the operator $C_F$ to
map continuously $A^p_{v,\alpha}({\Bbb D})$ into
$A^p_{u,\beta}(D)$ and the condition
$$N_{F,u,\beta}(z)=o(\gm_\al(\log|z|))\text{\qquad as\qquad}
|z|\to1$$ is necessary and sufficient for the operator
$C_F:\,A^p_{v,\alpha}({\Bbb D})\to A^p_{u,\beta}(D)$ to be
compact.
\eT
\begin{pf} Let us evaluate $N^*_{u,F,f,\beta}(w,r)$. We take the
set of all points $\{z_i\}$ in ${\Bbb D}$ such that $f(z_i)=w$ and
$v(z_i)>r$. Let $\cup_{j}A_{ij}=F^{-1}(z_i)$ be the decomposition
of the preimage of $z_i$ under $F$ into irreducible components.
The multiplicity of $f^*$ on $A_{ij}$ is equal to $m_im_{ij}$,
where $m_{ij}$ is the multiplicity of $F$ on $A_{ij}$ and $m_i$ is
the multiplicity of $f$ at $z_i$. If $r$ is so small that
$L(u_1)\cap T(r)=\emptyset$, where $T(r)=\{z\in D:\,v(F(z))>r\}$,
then by Proposition \ref{P:ln}
$$N^*_{u,F,f,\beta}(w,r)=
\sum_im_i\sum_jm_{ij}\int_{A_{ij}\cap
T(r)}\gm_\beta(u)(dd^cu)^{n-1}.$$ Considering $F$ as a holomorphic
function on $D_1$, we remark that
$$\sum_jm_{ij} \int_{A_{ij}\cap T(r)}\gm_\beta(u)(dd^cu)^{n-1}=
N_{F,\beta}(z_i),$$ so
$$N^*_{u,F,f,\beta}(w,r)=
\sum_im_iN_{F,\beta}(z_i).$$
\par We also observe that
$$N_{v,f,\al}(w)\ge \sum_im_i\gm_\al(v(z_i)).$$
Hence,
$$\frac{N^*_{u,F,f,\beta}(w,r)}{N_{v,f,\al}(w)}\le
\frac{\sum_im_iN_{F,\beta}(z_i)}{\sum_im_i\gm_\al(v(z_i))}\le
\max_i\left\{\frac{N_{F,\beta}(z_i)}{\gm_\al(v(z_i))}\right\}.$$
Thus,
$$\dl_{u,v,F,\beta,\al}(r)\le\sup_{|z|>r}
\frac{N_{F,\beta}(z)}{\gm_\al(v(z))}.$$ On the other hand if
$f(z)=z$, then
$$\frac{N^*_{u,F,f,\beta}(w,r)}{N_{v,f,\al}(w)}=
\frac{N_{F,\beta}(z)}{\gm_\al(v(z))}.$$ Hence,
$$\dl_{u,v,F,\beta,\al}(r)=\sup_{|z|>r}
\frac{N_{F,\beta}(z)}{\gm_\al(v(z))}.$$ Now  the sufficiency of
our conditions for $C_F$ to be bounded or compact  follows from
Theorem
\ref{T:bct}.
\par To show the necessity of the conditions we, firstly, note
that if $|F|<a<1$ on $D$, then the operator $C_F$ is both bounded
and compact and $N_{F,\beta}(z)=0$ when $|z|>a$. So this case is
trivial.
\par Now we assume that there is a sequence $\{\zeta_j\}\sbs D$
such that the points $z_j=F(\zeta_j)$ converge to the unit circle.
Following the standard argument (cf.  \cite{Sm}, \cite{CoM}) we
consider the functions
$$k_j(z)=\frac{(1-|z_j|^2)^{(\al+2)/p}}{(1-\ovr
z_jz)^{2(\al+2)/p}}$$ on ${\Bbb D}$. It can be shown (see
\cite{Sm}) that the norms $\|k_j\|_{A^p_\al}\approx1$ and the
functions $k_j$ converge to 0 uniformly on compacta.
\par Note that
$$k_j(z_j)=w_j=\frac{1}{(1-|z_j|^2)^{(\al+2)/p}}.$$
Thus,
$$
w_j^p\approx \frac{1}{\gamma_\alpha (\log |z_j|)}.
$$
Since $N_{u,k_j^\ast,\beta}(w_j)\ge N_{F,\beta}(z_j)$ this implies
$$
\nu_{F,\alpha,\beta}(w_j,k_j)\ge
C_1w_j^pN_{u,k_j^\ast,\beta}(w_j)\ge
C_2\frac{N_{F,\beta}(z_j)}{\gamma_\alpha (\log |z_j|)}.
$$
Since $k^\ast_j$ converges to 0 on any compact set $K\sbs D$, for a fixed
$a>1$ the condition $w_j>a \sup_{\zeta\in K}|k_j^\ast|(\zeta )$ is
satisfied and, therefore,
$$
\rho_{u_1,u_2,F,\alpha,\beta}(a)\ge
C_2\frac{N_{F,\beta}(z_j)}{\gamma_\alpha (\log |z_j|)}.
$$
For the same reason $\tilde{\rho}_{u_1,u_2,F,\alpha,\beta}(t,a)$
approaches 0 as $t\to 0$. Now the result follows  from Theorem
\ref{T:bctn}.
\end{pf}

\par We want to conclude this section with a useful
formula for the norms of composition operators when a hyperconvex
domain $D\sbs{\Bbb C}^n$ is mapped by a holomorphic function $F$
into the unit disk.
\bT\label{T:mud} Let $D\sbs{\Bbb C}^n$ be a hyperconvex domain
with an exhausting function $u$ such that the set $L(u)$ is
finite. If $F$ is a holomorphic function mapping $D$ into the unit
disk, $f$ is a holomorphic function on ${\Bbb D}$ and $f^*=C_Ff$,
then
$$\|f^*\|^p_{A^p_{u,\al}}=\il D\sigma_\al(u)|f|^p(dd^cu)^n+
\il{\Bbb D}N_{F,\al}(w)dd^c|f|^p.$$
\eT
\begin{pf} By Theorem \ref{T:cvi}
$$-\il{D}u\,dd^c|f^*|^p\wedge(dd^cu)^{n-1}=
\il{{\Bbb D}}N_{F,u}(w)\,dd^c|f|^p.$$ Thus, the case $\al=-1$
follows immediately from Theorem \ref{T:cvn}.
\par In the case $\al>-1$ we note that
\begin{equation}\begin{align}
&\il{-\infty}^0|r|^\al e^r\left(\il{B_u(r)}
(r-u)dd^c|f^*|^p\wedge(dd^cu)^{n-1}\right)\,dr\notag\\&=
\il{-\infty}^0|r|^\al e^r\il{{\Bbb
D}}N_{F,u}(w,r)\,dd^c|f|^p\,dr=\il{\Bbb D}N_{F,\al}(w)dd^c|f|^p.
\notag
\end{align}\end{equation}
Again, the identity follows from Theorem \ref{T:cvn}.
\end{pf}

\section{Composition operators induced by mappings into strongly
pseudoconvex domains}\label{S:comsp}

\par If $B$ is the unit ball in ${\Bbb C}^n$ with coordinates
$(z_1,\dots,z_n)$ and $u(z)=\log|z|$, then (see \cite[Prop.
1.4.10]{Ru1}) the function $\phi(z)=(1-z_1)^{-1}\in H^p_u(B)$ if
and only if $p<n$. So if $f:\,{\Bbb D}\to B$ is defined as
$f(z)=(z,0,\dots,0)$, then $C_f\phi\not\in H^p({\Bbb D})$ when
$p\ge 1$. On the other hand, the function $(1-z)^{-1}$ is in
$A^p_\al({\Bbb D})$ if and only if $p<\al+2$. So whatever is
$p<n$, the function $C_f\phi\in A^p_{n-2}({\Bbb D})$.
\par This calculation motivates the main result of this section.
Before proving it we mention that some special cases were
considered in \cite{MM}, \cite{CM}, \cite{KS} and \cite{SZ1}. A
similar theorem for mappings of polydisks was obtained in
\cite{SZ2}.
\bT\label{T:n-2}  Let $D_1\sbs{\Bbb C}^m $ be a hyperconvex domain and
$u\in\E(D_1)$. Let $D_2\sbs{\Bbb C}^m$ be a strongly pseudoconvex domain
with a $C^3$ exhausting strongly plurisubharmonic function $\rho$ such that
$\nabla\rho\ne0$ on $\bd D$. If $F: D_1\to D_2$ is a holomorphic mapping
then $C_F $ acts boundedly from $A^p_{\rho ,\alpha}(D_2)$ into
$A^p_{u,n+\alpha -1}(D_1), \ \alpha \geq -1$.
\eT

\begin{pf} Let $T_\zeta$ be the complex tangent plane to $\bd D_2$
at $\zeta\in\bd D_2$ and $U(\zeta,t)$ be the ball in $T_\zeta$ of
radius $\sqrt{t}$ centered at $\zeta$. Let
$$
A(\zeta,t)= \{ z\in D_2: \dist(z,U(\zeta,t)<t\}.
$$
A theorem of H\"ormander \cite{H} states  that if $\mu$ is a
positive measure on $D_2$ such that there is a constant $C$ so
that for every $\zeta\in \partial D$ and $t>0$,
\begin{equation}
\label{eq9} \mu(A(\zeta, t))<Ct^n,
\end{equation}
then there is a constant $C_1$ such that for every function $f\in
H^p_\rho(D_2)$
$$
\int_{D_2}|f|^pd\mu \leq C_1\| f\|_{H^p_\rho(D_2)}^p.
$$
The corresponding result for weighted Bergman spaces was proved by
Cima and Mercer (\cite{CM}). Namely, if
\begin{equation}
\label{eq9a} \mu (A(\zeta, t))<Ct^{\alpha +n+1},
\end{equation}
then there is a constant $C_1$ such that for every function $f\in
A^p_{\rho,\alpha}(D_2)$
$$
\int_{D_2}|f|^pd\mu \leq C_1\| f\|_{A^p_{\rho, \alpha}(D_2)}^p.
$$
\par By Corollary 3.2 from \cite[Ch. VII]{Ra} there are $\dl>0$
and a $C^2$ function $H(\zeta,w)$ on $\bd D\times D_\dl$, where
$D_\dl=\{\rho<\dl\}$, such that $H$ is holomorphic in $w$,
$|H(\zeta,w)|<1$ when $z\in\ovr D\sm\{\zeta\}$ and
$H(\zeta,\zeta)=1$. We may assume that the first and second
derivatives in $w$ of $H(\zeta,w)$ on $\bd D\times D_\dl$ do not
exceed some constant $C_1$. There is $t_0>0$ such that the sets
$U(\zeta,t)\sbs D_\dl$ when $0<t\le t_0$. Since the derivatives of
$H(\zeta,w)$ along the complex tangent plane $T_\zeta$ are equal
to 0 at $\zeta$, there is a constant $C_2$ such that
$|H(\zeta,w)-1|\le C_2t^{3/2}$ when $t\le t_0$ and $w\in
U(\zeta,t)$. Therefore, $|H(\zeta,w)-1|\le C_3t$ for some constant
$C_3$ when $t\le t_0$ and $w\in A(\zeta,t)$. It follows that there
is a constant $C_4$ such that the sets $A(\zeta,t)$ are contained
in the sets $F(\zeta,C_4t)=\{w\in D_2:\,|H(\zeta,w)-1|\le C_4t\}$.

\par Fix a point $z_0\in D_1$ and let $v(z)=g_{D_1}(z,z_0)$.
Write $h_\zeta(z)=H(\zeta,F(z))$ and consider the following
measure
$$\lm_\zeta(E)=
\il{-\infty}^0|r|^{n+\alpha -1}
e^r\mu_{v,r}(\chi_E(h_\zeta(z)))\,dr $$ on the unit disk ${\Bbb
D}$. If $\nu_F$ is a measure on $D_2$ defined by
$$\nu_F(E)=
\il{-\infty}^0|r|^{n+\alpha -1} e^r\mu_{v,r}(\chi_E(F(z)))\,dr,$$
then $\nu_F(A(\zeta,t))\le \lm_\zeta(E(1,C_4t))$ for all $t$.
Since there is $b<1$ such that $|h_\zeta(z_0)|\le b$ for all
$\zeta\in\bd D_2$, by Lemma
\ref{L:n1cm} there is a constant $a$ such that $\lm_\zeta(E(1,C_4t))\le
at^{n+\alpha +1}$. Thus, $\nu_F(A(\zeta,t))\le a_1t^{n+\alpha +1}$
for all $\zeta\in\bd D$ and all $t$.
\par Therefore, by theorems of H\"ormander and Cima and Mercer  we have for
$\phi\in A^p_{\rho,\alpha }(D_2)$
$$
\| C_F\phi\|_{A^p_{u,n+\alpha
-1}(D_1)}^p=\il{D_2}|\phi|^p\,d\nu_f\le C\|\phi\|_{A^p_{\rho,
\alpha}(D_2)}^p.
$$
\end{pf}

Under some additional assumptions about the mapping $F$ one might
get a considerable improvement of the result of the previous
Theorem. Below we consider two such cases.
\par In the first case $D\sbs{\Bbb C}^n$ is a strongly pseudoconvex domain
with a strictly plurisubharmonic exhaustion function $\rho\in C^3(V)$,
where $V$ is a neighborhood of $\ovr D$, and a holomorphic mapping $F$ is
defined on $V$ and takes $D$ into the unit disk ${\Bbb D}$.
\par We will need a lemma (see \cite[Ch. 3.13.37]{Sh}).
\bL\label{L:lc} There is $r_0>0$ and $C>0$ such that for every
point $z_0$ on the boundary of $D$ there are holomorphic
coordinates $(z_1,\dots,z_n)$ on $B(z_0,r_0)$ with the following
properties: $z_0=0$ and in these coordinates
\begin{equation}\label{e:ffr}
\rho(z)=\re z_n+\frac12H(z)+\phi(z),
\end{equation}
where
$$H(z)=\sum_{i,j=1}^n\frac{\bd^2\rho}{\bd z_j\bd\ovr
z_i}(0)z_i\ovr z_j$$ and $|\phi(z)|\le C|z|^3$.
\eL
\par We will call such coordinates {\it preferred}. In these
coordinates the real tangent hyperplane to $\bd D$ at $0$ is
$\{\re z_n=0\}$ and the complex tangent hyperplane to $\bd D$ at
$0$ is $\{z_n=0\}$.
\par The next lemma gives us a local estimate for the area of the
sets $\{F=w\}$.
\bL\label{L:lae} Suppose that $F$ is a holomorphic function on
a neighborhood $V$ of $\ovr D$ mapping $D$ into ${\Bbb D}$. Then there are
positive numbers $C$, $c$ and $r_1$ such that for every $z_0\in\bd D$ with
$F(z_0)=w_0$ and $|w_0|=1$ and for every $w\in{\Bbb C}$ and $r<0$ either
the set $\{F(z)=w\}\cap B(z_0,r_1)\cap\{\rho(z)<r\}$ is empty or
$r+c|w-w_0|\ge0$ and the area of this set does not exceed
$C(r+c|w-w_0|)^{(n-1)/2}$ when $r+c|w-w_0|<1$.
\eL
\begin{pf} We may assume that $w_0=1$, otherwise we consider
the function $\ovr w_0f(z)$.  Since $\rho(z)>0$ when $|F(z)|>1$ it
is easy to see that in preferred coordinates $\nabla
F(0)=(0,\dots,0,\lm)$, ${\im} \lm=0$ and $\lm>0$.
\par Let us write the Taylor expansion of $F$ as
$$F(z)=1+\lm z_n+a_0('z)+z_nb_0('z)+c_0z_n^2+F_3(z),$$
where $'z=(z_1,\dots,z_{n-1})$,
$$a_0('z)=\sum_{i,j=1}^{n-1}\frac{\bd^2F}{\bd z_j\bd
z_i}(0)z_i z_j,$$
$$b_0('z)=\sum_{j=1}^{n-1}\frac{\bd^2F}{\bd z_j\bd
z_n}(0)z_j,$$ $c_0=\bd^2F/\bd^2 z_n(0)$ and $|F_3(z)|<C_1|z|^3$
for some constant $C_1$.
\par By the Implicit Function Theorem there is a positive $r_1<r_0$
such that if $F(z)=1+w$ at some point $z\in B(0,r_1)$, then the
solution of the equation $F(z)=1+w$ in $B(0,r_1)$ can be
represented as $z_n=g('z,w)$, where $g$ is a holomorphic function.
We write $g$ as
\begin{equation}\label{eq:g}
g('z,w)=\lm^{-1}w+a('z)+wb('z)+cw^2+g_3('z,w)),\end{equation}
where
$$a('z)=\sum_{i,j=1}^{n-1}\frac{\bd^2g}{\bd z_j\bd
z_i}(0)z_i z_j,$$
$$b('z)=\sum_{j=1}^{n-1}\frac{\bd^2g}{\bd z_j\bd
w}(0)z_j,$$ $c=\bd^2g/\bd^2 w(0)$ and $|g_3('z,w)|<C_2|('z,w)|^3$
for some constant $C_2$.
\par Substituting (\ref{eq:g}) into the identity $F('z,g('z,w))=w$ we see that
$a('z)=-\lm^{-1}a_0('z)$, $b('z)=-\lm^{-2}b_0('z)$ and
$c=-\lm^{-3}c_0$. Thus
$$\rho_1('z,w)=\rho('z,g('z,w))=\lm^{-1}\re w+L('z,w)+\psi('z,w),$$
where
$$L('z,w)=-\re\left(\lm^{-1}a_0('z)+\lm^{-2}wb_0('z)+
\lm^{-3}c_0w^2\right)+\frac12H('z,\lm^{-1}w)$$ and
$|\psi'(z,w)|\le C_3|('z,w)|^3$ for some constant $C_3$.
\par Since $\rho_1('z,w)>0$ when $|1+w|^2>1$ or $2\re w>-|w|^2$,
we see that
$$L('z,w)\ge\frac1{2\lm}|w|^2.$$
\par Hence, the quadratic form $L$ is nonnegative and, therefore,
the set of vectors $('z,w)$, where $L('z,w)=0$ is a real linear
subspace. The quadratic form
$$L_1('z)=L('z,0)=-\frac1\lm\re a_0('z)+\frac12H('z,0)$$
is also nonnegative. Let us show that the real linear space $N$
where $L_1('z)=0$ does not contain complex lines. Observe that if
$L_1(\zeta 'z)=0$ for all $\zeta\in{\Bbb C}$, then
$$L_1('z)=-\frac1\lm\re \zeta a_0('z)+\frac12H(\zeta 'z,0).$$
Since $H$ is strictly positive we see that $L_1(\zeta 'z)>0$ when
$\zeta=-\ovr a_0('z)$. This contradiction shows that $N$ does not
contain complex lines and, therefore, its real dimension is at
most $n-1$.
\par The quadratic form
$L('z,w)=L_1('z)+\re wb_1('z)+c_1|w|^2$, where $b_1$ is a linear
function and $c_1>0$. Note that the coefficients of $b_1$ and
$c_1$ depend only on the second derivatives of $F$ and $\rho$ and,
therefore, are uniformly bounded. Hence if $\rho_1('z,w)<r$, then
there is a constant $c>0$ such that $L_1('z)\le r+c|w|$. Hence, if
the set $\{F(z)=w\}\cap B(z_0,r_1)\cap\{\rho(z)<r\}$ is non-empty,
then $r+c|w|\ge0$.
\par We can introduce orthonormal coordinates $x=(x_1,\dots,x_{2n-2})$ on
the ball in ${\Bbb C}^{n-1}$ such that
$$L_1(x)=\sum_{j=1}^kd_j|x_j|^2,$$
where all $d_j>0$ and $k\ge n-1$. The volume of the set $\{x\in
B:\,L_1(x)<r+c|w|\}$, where $B\sbs{\Bbb C}^{n-1}$ is the ball of
radius $r_1$ centered at the origin, does not exceed
$C_4(r+c|w|)^{k/2}$. Since the orthogonal projection of the set
$\{F=w\}$ in $B(0,r_1)$ has the Jacobian close to 1 we see that
the area of the set $\{F(z)=w\}\cap B(z_0,r_1)\cap\{\rho(z)<r\}$
does not exceed $C_4(r+c|w|)^{k/2}\le C_4(r+c|w|)^{(n-1)/2}$ when
$r+c|w|<1$.
\end{pf}
\par This lemma allows us to estimate the counting functions.
\bL\label{L:ecfs} In the assumptions of Lemma \ref{L:lae} there
are positive numbers $\dl$, $c>0$ and $C$ such that
$n_{F,\rho}(w,r)\le C(r+c|w|)^{(n-1)/2}$, $N_{F,\rho}(w,r)\le
C(r+c|w|)^{(n+1)/2}$ and $N_{F,\rho,\beta}(w)\le
C|w|^{\beta+(n+3)/2}$ when $|w|\ge 1-\eps$.
\eL
\begin{pf} Let us chose finitely many points $z_1,\dots,z_m$ in
the set $G=\{z\in\bd D:\,|F(z)|=1\}$ such that $G\sbs
W=\cup_{j=1}^mB(z_j,r_1)$. Clearly, there is $\eps>0$ such that
$z\in W$ when $w=f(z)$ and $|w|>1-\eps$. Since $\nabla F\ne0$ on
$G$ we may assume that $\eps$ is so small that the set
$\{f(z)=w\}$ is smooth when $|w|>1-\eps$. By Theorem \ref{T:cvi1}
$$n_{\rho,F}(w,r)=\il{\{f=w\}\cap B_\rho(r)}(dd^c\rho)^{n-1}.$$
But $dd^c\rho$ is equivalent to the Euclidean metric and, there
fore, $n_{\rho,F}(w,r)$ does not exceed the area of the set
$\{f(z)=w\}\cap B_\rho(r)$ times some constant. In every ball
$B(z_{j},r_1)$ this area does not exceed $C(r+c|w|)^{(n-1)/2}$.
Hence $n_{\rho,F}(w,r)\le mC(r+c|w|)^{(n-1)/2}$.
\par The function
$$N_{\rho,F}(w,r)=\il{-c|w|}^rn_{\rho,F}(w,t)\,dt\le
\frac{2mC}{n+1}(r+c|w|)^{(n+1)/2}.$$
\par Finally,
$$N_{\rho,F,\beta}(w)=\il{-c|w|}^0|r|^\beta e^rN_{\rho,F}(w,t)\,dt\le
mC\il{-c|w|}^0|r|^\beta(r+c|w|)^{(n+1)/2}\,dt.$$ To estimate the
last integral we notice that
$$\il{-c|w|}^{-c|w|/2}|r|^\beta(r+c|w|)^{(n+1)/2}\,dt\le
C_1|w|^{\beta+(n+3)/2}$$ and
$$\il{-c|w|/2}^{0}|r|^\beta(r+c|w|)^{(n+1)/2}\,dt\le
C_2|w|^{\beta+(n+3)/2}.$$ Thus $N_{\rho,F,\beta}(w)\le
C_3|w|^{\beta+(n+3)/2}$.
\end{pf}
\par The following theorem is an immediate consequence of these
lemmas and Theorem \ref{T:bctd}.
\bT\label{T:spdd}  Suppose that $F$ is a holomorphic function on
a neighborhood $V$ of $\ovr D$ mapping a strongly pseudoconvex
domain $D$ into ${\Bbb D}$. Then the composition operator $C_F$
maps $A^p_\al({\Bbb D})$ into $A^p_\beta(D)$ when
$\al\le\beta+(n-1)/2$. This operator is compact when
$\al<\beta+(n-1)/2$.
\eT
\par It is interesting to note that this theorem cannot be
improved even for quadratic polynomials. Let $B$ be the unit ball
in ${\Bbb C}^n$ centered at the origin and
$F(z)=z_1^2+\dots+z_n^2$. If $z_j=r_je^{i\phi_j}$ then $|F(z)|=1$
if and only if $\sum r_j^2=1$ and $e^{2i\phi_j}=1$. Let us
estimate the area of the set $\{F=w\}$ near $z=(0,\dots,0,1)$. The
equation $F(z)=1+w$ has a holomorphic solution $g('z,w)$ with the
following Taylor expansion of order 2
$$g('z,w)=1+\frac w2-\frac{'z^2}2-\frac{w^2}4,$$ where
$'z^2=z_1^2+\dots+z_{n-1}^2$.
\par If $|'z|^2+|g('z,w)|^2<1+r$, $r<0$, then
$\re w-\re 'z^2+|'z|^2<cr$ for some constant $c>0$. If
$z_j=x_j+iy_j$ this inequality is equivalent to
$y_1^2+\dots+y_{n-1}^2\le cr-\re w$. Now the estimates similar to
what was used in the proof of Lemma \ref{L:lae} tell us that
$n_{\rho,F}(w,r)\ge C(cr+|w|)^{(n-1)/2}$. Consequently,
$N_{\rho,F,\beta}(w)\ge C|w|^{\beta+(n+3)/2}$.

\vspace{.1cm}

\par  Our next example shows that Hardy and Bergman spaces
are preserved by proper mappings of domains with equal dimension.
\par Suppose that domains $D_1$ and $D_2$ are in ${\Bbb C}^n$ and
$f:\,D_1\to D_2$ is a proper holomorphic mapping, i.e.,
$f^{-1}(K)$ is compact when $K\sbs D_2$ is compact. In particular,
for every $w\in D_2$ the set $\{f^{-1}(w)\}$ is compact and
analytic, so it is finite. The branch set $B_f$ of $f$ is the set
of points in $D_1$, where $f$ is not locally homeomorphic. It is
contained in the analytic set $\{J_f=0\}$, where $J_f$ is the
Jacobian of $f$ and, consequently, has the dimension less than
$n$. Thus, if points $w_1,w_2\in D_2\sm f(B_f)$ then we can
connect them by a continuous curve $\gm\subset D_2\setminus
f(B_f)$ and considering its preimage to see that the sets
$\{f^{-1}(w_1)\}$ and $\{f^{-1}(w_2)\}$ have the same number of
points $m$. This number is called the {\it multiplicity} of $f$.
\par If $z_0\in D_1$ and $w_0=f(z_0)$ we consider a ball
$B_r\sbs\sbs D_2$ centered at $w_0$ and of radius $r$. If $W$ is a
connected component of $f^{-1}(B_r)$ containing $z_0$, then the
restriction of $f$ to $W$ maps $W$ properly on $B_r$ with the
multiplicity $m_r$. The limit of $m_r$ as $r$ goes to 0 is called
the {\it multiplicity} of $f$ at $z_0$ and denoted by $m(z_0,f)$.
\par If $\phi$ is a function on $D_2$, then as above we set
$\phi^*(z)=\phi(f(z))$. If $\phi$ is a function on $D_1$ then we
let
$$\phi_*(w)=\frac1m\sum_{f(z)=w}\phi(z),$$
where the summation counts multiplicities. Clearly, $\phi_*$ is
plurisubharmonic when $\phi$ is plurisubharmonic and holomorphic
when $\phi$ is holomorphic.
\par The following theorem is well-known.
Unfortunately, we could not find a reference so we give a proof.
\bT\label{T:pmi} Suppose that domains $D_1$ and $D_2$ are in ${\Bbb C}^n$ and
$f:\,D_1\to D_2$ is a proper holomorphic mapping of multiplicity
$m$.
\par If $v_1,\dots,v_n$ are plurisubharmonic functions on $D_2$
and $\phi$ is a non-negative Borel function on $D_2$, then
$$m\il{D_2}\phi dd^cv_1\wedge\dots\wedge dd^cv_n=
\il{D_1}\phi^* dd^cv^*_1\wedge\dots\wedge dd^cv^*_n.$$
\par If $\phi$ is a non-negative Borel function on $D_1$, then
$$m\il{D_2}\phi_* dd^cv_1\wedge\dots\wedge dd^cv_n=
\il{D_1}\phi^* dd^cv^*_1\wedge\dots\wedge dd^cv^*_n.$$
\eT
\begin{pf} To prove the first part we take an exhaustion of $D_2$
by domains $V_n$ such that $V_k\sbs\sbs V_{k+1}\sbs\sbs D_2$. It
is known that for each $k$ there are smooth plurisubharmonic
functions $v_{jk}$, $1\le j\le n$, defined on $V_{k+1}$ and such
that $v_{jk}\ge v_{j,k+1}$ on $V_{k}$ and $\lim v_{jk}=v_j$. We
also take continuous non-negative functions $h_k\le1$ on $D_2$
such that $h_k\equiv1$ on $V_k$ and $h_k\equiv0$ on $D\sm
V_{k+1}$.

\par By the result of Bedford and Taylor (\cite{Kl}, p.114)
$$\lim_{k\to\infty}\il{D_2}h_l\phi dd^cv_{1k}\wedge\dots\wedge dd^cv_{nk}=
\il{D_2}h_l\phi dd^cv_1\wedge\dots\wedge dd^cv_n.$$
\par Since the sets $f(B_f)$ and $f^{-1}(f(B_f))$ have a zero
measure,
$$m\il{D_2}h_l\phi dd^cv_{1k}\wedge\dots\wedge dd^cv_{nk}=
\il{D_1}h_l\phi dd^cv_{1k}^*\wedge\dots\wedge dd^cv_{nk}^*.$$
Taking first the limit as $k\to\infty$ and then the limit as
$l\to\infty$ we get our statement.
\par The second statement has a similar proof.
\end{pf}
\par Now we can prove a theorem about the composition operators of
proper holomorphic mappings.
\bT\label{T:pm} Let $D_1$ and $D_2$ be hyperconvex domains with
exhaustion functions $u_j\in\E(D_j)$, $j=1,2$, and let $f:\,D_1\to
D_2$ be a proper holomorphic mapping. Then the composition
operator $C_f$ maps $PS_{u_2}(D_2)$ into $PS_{u_1}(D_1)$ and
$A_{u_2,\al}(D_2)$ into $A_{u_1,\al}(D_1)$. Moreover, there is a
constant $A\ge1$ such that
$A^{-1}\|\phi\|_{u_2}\le\|C_f\phi\|_{u_1}\le A\|\phi\|_{u_2}$ and
$A^{-1}\|\phi\|_{u_2,\al}\le\|C_f\phi\|_{u_1,\al}\le
A\|\phi\|_{u_2,\al}$.
\par Consequently, $C_f$ maps continuously Banach spaces $H^p_{u_2}(D_2)$
into $H^p_{u_1}(D_1)$ and $A^p_{u_2,\al}(D_2)$ into
$A^p_{u_1,\al}(D_1)$.
\eT
\begin{pf} We fix a point $w_0\in D_2$. Then we consider
the exhaustion functions $v(z)=g_{D_2}(w,w_0)$ on $D_2$ and
$v^*(z)$ on $D_1$. Both functions belong to $\E(D_2)$ and
$\E(D_1)$ respectively.
\par Note that $(v^*)_r=(v_r)^*$ and $(dd^c v)^n\equiv0$ outside
of $\ovr B_v(r)$ while $(dd^c v^*)^n\equiv0$ outside of $\ovr
B_{v^*}(r)$. Thus
$$\mu_{v,r}(\phi)=\il{D_2}\phi(dd^cv_r)^n\text{ and }
\mu_{v^*,r}(\phi^*)=\il{D_1}\phi^*(dd^cv^*_r)^n.$$ By Theorem
\ref{T:pmi} $m\mu_{v,r}(\phi)=\mu_{v^*,r}(\phi^*)$, where $m$ is
the multiplicity of $f$. Hence, if $\phi$ is in $PS_v(D_2)$ or in
$A_{v,\al}(D_2)$, then $\phi^*$ in $PS_{v^*}(D_1)$ or in
$A_{v^*,\al}(D_1)$ and $\|\phi^*\|_{v^*}=m\|\phi\|_v$ and
$\|\phi^*\|_{v^*,\al}=m\|\phi\|_{v,\al}$ respectively.
\par Now our theorem follows immediately from Theorem \ref{T:ni}.
\end{pf}
\par The operator $L_f$ mapping the functions $\phi$ on $D_1$ into
the functions $\phi_*$ on $D_2$ is linear and $L_f\circ C_f$ is
the identity operator. As the next theorem shows this operator has
the same properties as $C_f$ and the proof is also the same.
\bT\label{T:pm2} In the assumptions of Theorem \ref{T:pm} the
operator $L_f$ maps $PS_{u_1}(D_1)$ into $PS_{u_2}(D_2)$ and
$A_{u_1,\al}(D_1)$ into $A_{u_2,\al}(D_2)$. Moreover, there is a
constant $A\ge1$ such that
$A^{-1}\|\phi\|_{u_1}\le\|L_f\phi\|_{u_2}\le A\|\phi\|_{u_1}$ and
$A^{-1}\|\phi\|_{u_1,\al}\le\|L_f\phi\|_{u_12,\al}\le
A\|\phi\|_{u_1,\al}$.
\par Consequently, $L_f$ maps continuously Banach spaces $H^p_{u_1}(D_1)$
into $H^p_{u_2}(D_2)$ and $A^p_{u_1,\al}(D_1)$ into
$A^p_{u_2,\al}(D_2)$.
\eT

\end{document}